\documentclass[11pt]{article}
\usepackage[dvips]{graphics}
\input{amssym.tex}

\usepackage{xcolor}

\title{\Large\bf 
DYNAMIC PROBLEM OF A POWER-LAW GRADED HALF-PLANE AND AN ASSOCIATED CARLEMAN  PROBLEM FOR TWO FUNCTIONS}
\author{Y.A. Antipov\\
Department of Mathematics, Louisiana State University,\\
Baton Rouge LA 70803, USA\\
}

\newcommand{\sgn}{\mathop{\rm sgn}\nolimits}

\newcommand{\R}{\mathop{\rm Re}\nolimits}
\newcommand{\const}{\mbox{const}}

\newcommand{\beqa}{\begin{eqnarray}}
\newcommand{\eeqa}[1]{\label{#1}\end{eqnarray}}
\newcommand{\bequ}{\begin{equation}}
\newcommand{\eequ}[1]{\label{#1}\end{equation}}

\newcommand{\Md}{\partial}

\newcommand{\ov}[1]{\overline{#1}}

\newcommand{\Ga}{\alpha}
\newcommand{\Gb}{\beta}
\newcommand{\Gd}{\delta}

\newcommand{\Gve}{\varepsilon}
\newcommand{\Gf}{\phi}
\newcommand{\Gvf}{\varphi}
\newcommand{\Gg}{\gamma}

\newcommand{\Gk}{\kappa}
\newcommand{\Gvk}{\varkappa}
\newcommand{\Gl}{\lambda}
\newcommand{\Gn}{\eta}
\newcommand{\Gm}{\mu}

\newcommand{\Gt}{\theta}

\newcommand{\Gr}{\rho}

\newcommand{\Gs}{\sigma}

\newcommand{\Gx}{\xi}

\newcommand{\GD}{\Delta}
\newcommand{\GF}{\Phi}
\newcommand{\GG}{\Gamma}

\newcommand{\GP}{\Pi}

\newcommand{\GO}{\Omega}

\newcommand{\CF}{{\cal F}}

\newcommand{\CI}{{\cal I}}
\newcommand{\CJ}{{\cal J}}

\newcommand{\CM}{{\cal M}}

\def\BF{{\bf F}}

\newcommand{\beq}{\begin{equation}}
\newcommand{\eeq}{\end{equation}}
\newcommand{\barr}{\begin{eqnarray}}
\newcommand{\earr}{\end{eqnarray}}
\newcommand{\beqn}{\begin{equation*}}
\newcommand{\eeqn}{\end{equation*}}
\newcommand{\barrn}{\begin{eqnarray*}}
\newcommand{\earrn}{\end{eqnarray*}}

\newcommand{\fr}{\frac}

\begin{document}
\maketitle

\begin{abstract}

A steady state plane problem of an inhomogeneous half-plane subjected to a load running along the boundary at subsonic speed is analyzed. The Lame coefficients and the density of the half-plane are assumed to be power functions of depth. The model is different from the standard static model have been used in contact mechanics since the Sixties and originated from the 1964 Rostovtsev exact solution of the Flamant problem of a power-law graded half-plane. To solve the governing dynamic equations with variable coefficients written in terms of the displacements, we propose a method that, by means of the Fourier and Mellin transforms, maps the model problem to a Carleman boundary value problem for two meromorphic functions in a strip with two shifts or, equivalently, to a system of two difference equations of the second order with variable coefficients.  By partial factorization the Carleman  problem is recast as a system of four singular integral equations on a segment with a fixed singularity and highly oscillating coefficients. A numerical method for its solution is proposed and tested. Numerical results for the displacement and stress fields are presented and discussed. 

\end{abstract}

\setcounter{equation}{0}

\section{Introduction} 

The standard model used in contact mechanics to describe the indentation of a stamp into a power-law graded 
elastic body is due to Lekhnitskii {\bf(\ref{lek})} and Rostovtsev  {\bf(\ref{ros})}. The former paper uses the method of separation of variables to derive a solution
of the two-dimensional problem
 of finding the distribution of the material constants which admits the given state of stress in a wedge subjected to a concentrated force applied at the wedge vertex. Specifically, it points out
 that, when the Poisson ratio $\nu_P$ is constant, the Young modulus is $E(r)=E_0 r^\nu$, and $E_0$ is constant,  the stresses are distributed radially.
An exact solution of the Flamant problem of a half-plane in the same case, $\nu_P=\const$ and  $E(r)=E_0 r^\nu$, is derived in   {\bf(\ref{ros})}.
The results include the integral representation of the normal displacement on the boundary of the half-plane 
\beq
v(x)=\Gt_0\int_a^b\fr{p(\Gx)d\Gx}{|x-\Gx|^\nu},
\label{1.1}
\eeq
where  $p(\Gx)$ is the pressure distribution, $(a,b)$ is the contact zone, and the parameter $\Gt_0$ is expressed through $\nu$, $\nu_P$, $E_0$, and the magnitude  of the normal force in the Flamant model. 
Equation (\ref{1.1}) and its  axisymmetric analog  {\bf(\ref{kor})},  {\bf(\ref{mos})}, {\bf(\ref{pop1})} had been used  and solved exactly even before the papers    {\bf(\ref{lek})} and {\bf(\ref{ros})} were published. In those preliminary studies the parameter
$\Gt_0$ was not specified explicitly and proposed to be determined experimentally. 
The Rostovtsev solution   {\bf(\ref{ros})} was employed to  model the plane and axisymmetric Hertzian  contact interaction of a rigid stamp
with a power-graded foundation {\bf(\ref{pop2})}. Later, in concert with the  Johnson-Kendall-Roberts 
model {\bf(\ref{joh})}, the  Rostovtsev model was used in  {\bf(\ref{gia})} -- {\bf(\ref{hes})}.
By employing the same model the Hertzian and Johnson-Kendall-Roberts  contact interaction of two  elastic power-law graded semi-infinite bodies   $\{z<0\}$ and $\{z>0\}$
characterized by the Young moduli $E(z)=E_1 |z|^{\nu_1}$ and $E(z)=E_2 z^{\nu_2}$ was studied in the 
 axisymmetric case  in {\bf(\ref{ant1})} and  in the plane case in {\bf(\ref{ant2})}. Recently, in the framework of the Rostovtsev model,
  the plane and axisymmetric fracture problems of an interfacial finite crack between  two power-law graded materials were solved 
  exactly {\bf(\ref{ant3})} by the method of orthogonal polynomials and the Wiener-Hopf method.   
 
 All these models are originated from the Rostovtsev solution  {\bf(\ref{ros})} of the Flamant problem of a power-law graded half-plane  
and share two shortcomings. Firstly, they cannot describe the stress distribution and displacements outside the contact zone.
Secondly, the models are static and do not admit a generalization to the dynamic case.  The main goal of the present work is to 
design an {\it ab initio} model capable to recover the displacement and stress fields everywhere in the elastic body 
in both static and dynamic cases. To do this, we aim to solve the dynamic boundary value problem of a power-law graded
half-plane subject to loading running on the boundary at constant speed. We confine ourselves to considering the steady state subsonic regime.

In Section 2 we write down the dynamic boundary value problem  with respect to the tangential and normal displacements in a half-plane when the two Lame constants and the density are functions of depth $y$, 
$\Gl(y)=\Gl_0 y^\nu$,  $\Gm(y)=\Gm_0 y^\nu$, and   $\Gr(y)=\Gr_0 y^\nu$,  $0<y<\infty$.
This results in variation with $y$ of  the coefficients of the governing system of partial differential equations written in terms of the displacements.
On the boundary, we assume that the two traction components $\Gs_{j2}$, $j=1,2$,
are prescribed to be $h_{j}(\Gx)$, $\Gx=x-Vt$, $V$ is speed, or, equivalently, $y^\nu\Md u_1/\Md y\to \mu_0^{-1}h_1(\Gx)$ and 
$y^\nu\Md u_2/\Md y\to (\Gl_0+2\mu_0)^{-1}h_2(\Gx)$, $y\to 0^+$. 

In Section 3 we apply the Fourier and Mellin transforms 
to map the boundary value problem to  a system of two difference equations of the second order with variable 
coefficients or, equivalently, to the Carleman boundary value problem 
{\bf(\ref{car})},  {\bf(\ref{zve})} with two shifts in a strip for two meromorphic functions.  In the scalar case and when there is only one shift, the problem in a strip is equivalent to a scalar Riemann-Hilbert problem and admits a closed-form solution {\bf(\ref{che})}, {\bf(\ref{ban})}.  The two-shift-scalar-problem with periodic 
coefficients and the one-shift-problem for two functions with a special matrix coefficient also admit closed-form solutions  {\bf(\ref{ant4})},
  {\bf(\ref{ant5})} by the method of Riemann surfaces. The problem derived in Section 3 cannot be solved in closed form 
  by the methods available in the literature. 
By applying the method of partial factorization we recast  the problem as a system of four singular integral equations with oscillating coefficients and arbitrary constants in the right-hand side.
These constants are fixed by satisfying the conditions on the boundary of the elastic half-plane  in Section 4. 

In Section 5 we analyze the particular case when there are no body forces and  the loading applied is a concentrated 
force $\Gs_{j2}(\Gx,0)=H_j\Gd(\Gx)$, $|\Gx|<\infty$, $\Gd(\Gx)$ is the Dirac function. We also derive integral representations and asymptotic expansions of the 
displacements and stresses. A  method for numerical solution of the system of integral equations is proposed in Section 6. 
Notice that similar integral equations arise in diffraction theory  {\bf(\ref{ant6})},  {\bf(\ref{ant7})}. However, there is a significant
difference between those equations and the ones solved in the present paper. In the case of the dynamic elastic problem of a power-law
graded    half-plane the coefficients of the integral equations oscillate at one of the  endpoints and do not have definite limits, while
in the diffraction problems such limits exist. 
Numerical results for the displacements and stresses are presented and discussed in Section 7.
 
\setcounter{equation}{0}

\section{Formulation}\label{form}

An inhomogeneous elastic half-plane $\{|x|<\infty, y>0\}$ is considered, with its boundary subjected to loading $h_j(x-Vt)$
running along the boundary at constant speed $V$
\beq
\Gs_{j2}(x,0,t)=h_j(x-Vt), \quad |x-Vt|<\infty, \quad j=1,2.
\label{2.1}
\eeq
 Here, $\Gs_{j2}$ are the traction components and $h_j$  are prescribed functions. The Lame coefficients $\Gl$ and $\Gm$ and the mass density $\Gr$
are independent of $x$ and time $t$ and are power functions of depth,
\beq
\Gl(y)=\Gl_0 y^\nu, \quad  \Gm(y)=\Gm_0 y^\nu, \quad  \Gr(y)=\Gr_0 y^\nu, \quad y>0,
\label{2.2}
\eeq
where $\Gl_0$, $\Gm_0$, and $\Gr_0$ are positive constants, and  $0<\nu<1$. 
The momentum balance equations of two-dimensional dynamic elasticity have the form 
$$
\Gs_{11,1}+\Gs_{12,2}+\Gr f_1=\Gr \ddot u_1,
$$
\beq
\Gs_{12,1}+\Gs_{22,2}+\Gr f_2=\Gr \ddot u_2, \quad |x|<\infty, \quad y>0, \quad t>0,
\label{2.3}
\eeq
and the stress-strain relations written in terms of the stresses and the displacements derivatives are
$$
\Gs_{11}=(\Gl+2\mu)u_{1,1}+\Gl u_{2,2},
\quad 
\Gs_{22}=\Gl u_{1,1}+(\Gl+2\mu)u_{2,2},
$$
\beq
\Gs_{12}=\mu(u_{1,2}+u_{2,1}), \quad |x|<\infty, \quad y\ge 0, \quad t>0.
\label{2.4}
\eeq 
Here, $f_j=f_j(x-Vt,y)$ are body forces per unit mass,  the notations $g_{,1}$  and $g_{,2}$ mean the partial derivatives $g_x$ and $g_y$, respectively, while $\ddot u_j$ denotes the second time-derivative of  $u_j$.

To proceed with the solution, we eliminate the stresses derivatives from the momentum balance equations and obtain the following equations governing the dynamics of a power-law graded material:
$$
(\Gl_0+2\mu_0)u_{1,11}+\mu_0 u_{1,22}+(\Gl_0+\mu_0)u_{2,12}+\fr{\mu_0\nu}{y}(u_{1,2}+u_{2,1})+\Gr_0f_1=\Gr_0 \ddot u_1,
$$    
\beq
\mu_0 u_{2,11}+(\Gl_0+2\mu_0)u_{2,22}+(\Gl_0+\mu_0)u_{1,12}+\fr{\nu}{y}[\Gl_0u_{1,1}+(\Gl_0+2\mu_0)u_{2,2}]+\Gr_0f_2=\Gr_0 \ddot u_2.
\label{2.5}
\eeq
Since the time-dependence of the loading and body force functions is realized through the variable $\Gx=x-Vt$, the
problem is steady state, and the mechanical fields are functions of two variables, $\Gx$ and $y$,
\beq
\Gs_{ij}=\Gs_{ij}(\Gx,y), \quad u_j=u_j(\Gx,y), \quad i,j=1,2.
\label{2.6}
\eeq
In the steady state conditions the displacements are found to satisfy the equations
 $$
(\Gl_0+2\mu_0-\Gr_0 V^2)\fr{\Md^2u_1}{\Md\Gx^2}+\mu_0 \fr{\Md^2 u_1}{\Md y^2}+(\Gl_0+\mu_0)\fr{\Md^2 u_2}
{\Md\Gx\Md y}+\fr{\mu_0\nu}{y}\left(\fr{\Md u_1}{\Md y}+\fr{\Md u_2}{\Md\Gx}\right)+\Gr_0 f_1=0,
$$    
$$
(\mu_0-\Gr_0 V^2)\fr{\Md^2u_2}{\Md\Gx^2}+(\Gl_0+2\mu_0) \fr{\Md^2 u_2}{\Md y^2}+(\Gl_0+\mu_0)\fr{\Md^2 u_1}
{\Md\Gx\Md y}+\fr{\nu}{y}\left[\Gl_0\fr{\Md u_1}{\Md\Gx}+(\Gl_0+2\Gm_0)\fr{\Md u_2}{\Md y}\right]+\Gr_0 f_2=0.
$$
\beq
-\infty<\Gx<\infty, \quad 0<y<\infty.
\label{2.7}
\eeq
To complete the formulation, we need to express the boundary conditions in terms of the displacements derivatives.
From the relations (\ref{2.2}) and (\ref{2.4}) we have
\beq
\Gs_{12}=\mu_0 y^\nu\left(\fr{\Md u_1}{\Md y}+\fr{\Md u_2}{\Md\Gx}\right),
\quad 
\Gs_{22}=y^\nu\left[\Gl_0\fr{\Md u_1}{\Md \Gx}+(\Gl_0+2\mu_0)\fr{\Md u_2}{\Md y}\right].
\label{2.8}
\eeq
In the classical theory, when the Lame parameters are constants, both of the derivatives of the displacements, the tangential and the 
normal ones, contribute to the boundary conditions on the boundary of the half-plane. For a power-law graded half-plane, because of the term $y^\nu$ with $0<\nu<1$,
and because the tangential derivatives of the displacements are bounded as $y\to 0^+$,
the boundary conditions read
\beq
\mu_0\lim_{y\to 0^+} y^\nu\fr{\Md u_1}{\Md y}=h_1(\Gx), \quad (\Gl_0+2\mu_0)\lim_{y\to 0^+} y^\nu\fr{\Md u_2}{\Md y}=h_2(\Gx),
\quad |\Gx|<\infty.
\label{2.9}
\eeq
Equations (\ref{2.7}) and the conditions (\ref{2.9}) constitute the  boundary value problem of the model to be solved.

\setcounter{equation}{0}

\section{Carleman boundary value problem for two meromorphic  functions in a strip}\label{car}

Applying the integral Fourier transform with respect to $\Gx$
\beq
\tilde u_j(\Ga,y)=\int_{-\infty}^\infty u_j(\Gx,y)e^{i\Ga\Gx}d\Gx, \quad 
\tilde f_j(\Ga,y)=\int_{-\infty}^\infty f_j(\Gx,y)e^{i\Ga\Gx}d\Gx,\quad j=1,2,
\label{3.1}
\eeq
reduces the dimension of the problem. We have
$$
\mu_0\left(\fr{d^2 \tilde u_1}{dy^2}+\fr{\nu}{y}\fr{d \tilde u_1}{dy}\right)-(\Gl_0+2\mu_0-\Gr_0 V^2)\Ga^2\tilde u_1-(\Gl_0+\Gm_0)i\Ga\fr{d \tilde u_2}{dy}-
\fr{\mu_0 \nu i\Ga}{y} \tilde u_2+\Gr_0\tilde f_1=0,
$$
$$
(\Gl_0+2\mu_0)\left( \fr{d^2 \tilde u_2}{dy^2}+\fr{\nu}{y}\fr{d \tilde u_2}{dy}\right)-(\mu_0-\Gr_0 V^2)\Ga^2\tilde u_2-(\Gl_0+\mu_0)i\Ga\fr{d \tilde u_1}{dy}-
\fr{\Gl_0 \nu i\Ga}{y} \tilde u_1+\Gr_0\tilde f_2=0,
$$
\beq
0<y<\infty.
\label{3.2}
\eeq
It follows from (\ref{2.9}) that the Fourier-transforms of the displacements $\tilde u_j(\Ga,y)$ have to satisfy the boundary conditions
\beq
\mu_0\lim_{y\to 0^+} y^\nu\fr{d \tilde u_1}{d y}=\tilde h_1(\Ga), \quad (\Gl_0+2\mu_0)\lim_{y\to 0^+} y^\nu\fr{d \tilde u_2}{d y}=\tilde h_2(\Ga),
\label{3.3}
\eeq
where we denoted 
\beq
\tilde h_j(\Ga)=\int_{-\infty}^\infty h_j(\Gx)e^{i\Ga\Gx}d\Gx,\quad j=1,2.
\label{3.4}
\eeq
Next we intend to apply the Mellin transform with respect to $y$
\beq
\hat u_j(\Ga,s)=\int_0^\infty \tilde u_j(\Ga,y)y^{s-1}dy.
\label{3.5}
\eeq
We seek the displacements in the class of functions bounded as $y\to 0^+$ and vanishing at infinity as $y^{-\Gb}$ ($0<\Gb<1$), that is
\beq
\tilde u_j(\Ga,y)=O(1), \quad y\to 0^+, \quad \tilde u_j(\Ga,y)=O(y^{-\Gb}), \quad y\to\infty.
\label{3.6}
\eeq
This guarantees that the functions $\hat u_j(\Ga,s)$ are holomorphic in the strip $0<\R s<\Gb$.
Continue analytically these functions to meromorphic functions in a wider strip $\GP=\{\Gs-2<\R s<\Gs\}$, where
$\Gs\in(0,\Gb)$. It immediately follows from the boundedness of the displacements and the boundary conditions (\ref{3.3})
that
\beq
\tilde u_j(\Ga,y)=B'_j(\Ga) +B_j''(\Ga)y^{1-\nu}+B_j(\Ga,y), \quad y\to 0^+,
\label{3.7}
\eeq
with $B_j'$ and $B_j''$ being independent of $y$, while $B_j(\Ga,y)$ being such that $B_j(\Ga,y)=o(y^{1-\nu})$, $y\to 0^+$. We split the Mellin integrals (\ref{3.5}) into integrals
over $(0,1)$ and $(1,\infty)$,  assume  that $\R s>0$ and  integrate by parts the finite integrals. We have 
\beq
\hat u_j(\Ga,s)=\fr{B'_j(\Ga)}{s}+\fr{B_j''(\Ga)}{s-\nu+1}+\int_0^1 B_j(\Ga,y)y^{s-1}dy+\int_1^\infty \tilde u_j(\Ga,y)y^{s-1}dy.
\label{3.8}
\eeq
The functions $\hat u_j(\Ga,s)$ are meromorphic functions  with respect to $s$ in the strip $\GP$.
They have simple poles at the points $s=0$ and $s=\nu-1$. 
To apply the Mellin transforms to the derivatives $d \tilde u_j/dy$ and $d^2 \tilde u_j/dy^2$ we assume 
first that $\R s> 1$ and $\R s>2$, respectively, integrate by parts and then extend the result analytically to the contour
$\GO=\{\R s=\Gs\in (0,\Gb)\}$,
$$
\int_0^\infty \fr{d\tilde u_j(\Ga,y)}{dy} y^{s-1}dy=-(s-1)\hat u_j(\Ga,s-1), \quad s\in\GO,
$$
\beq
\int_0^\infty \fr{d^2\tilde u_j(\Ga,y)}{dy^2} y^{s-1}dy=(s-1)(s-2)\hat u_j(\Ga,s-2), \quad s\in\GO.
\label{3.9}
\eeq
We are now ready to write down the Mellin images of equations (\ref{3.2}). By using formulas (\ref{3.9}) we deduce
$$
\Ga^2(\Gr_0 V^2-\Gl_0-2\mu_0)\hat u_1(\Ga,s)+i\Ga [(\Gl_0+\mu_0)(s-1)-\mu_0\nu] \hat u_2(\Ga,s-1)
$$
$$
+\mu_0(s-2)(s-1-\nu)\hat u_1(\Ga,s -2)=-\Gr_0\hat f_1(\Ga,s), \quad s\in\GO.
$$
$$
\Ga^2(\Gr_0 V^2-\mu_0)\hat u_2(\Ga,s)+i\Ga [(\Gl_0+\mu_0)(s-1)-\Gl_0\nu] \hat u_1(\Ga,s-1)
$$
\beq
+(\Gl_0+2\mu_0)(s-2)(s-1-\nu)\hat u_2(\Ga,s -2)=-\Gr_0\hat f_2(\Ga,s), \quad s\in\GO.
\label{3.10}
\eeq
Introduce the shear and wave speeds
\beq
c_d=\sqrt{\fr{\Gl+2\mu}{\Gr}}=\sqrt{\fr{\Gl_0+2\mu_0}{\Gr_0}}, \quad 
c_s=\sqrt{\fr{\mu}{\Gr}}=\sqrt{\fr{\mu_0}{\Gr_0}},
\label{3.11}
\eeq
where
\beq
\Gl_0=\fr{E_0\nu_P}{(1+\nu_P)(1-2\nu_P)}, \quad \mu_0=\fr{E_0}{2(1+\nu_P)},
\label{3.12'}
\eeq
$\nu_P$ is the Poisson ratio, $E=E_0 y^{\nu}$ is the Young modulus, We also denote
\beq
a_d=\fr{c_d}{V}, \quad a_s=\fr{c_s}{V}.
 \label{3.12}
 \eeq
 In these notations, the equations can be rewritten as a system of two difference equations of the second order with variable coefficients
 or, equivalently, as the following Carleman boundary value problem for a strip $\GP=\{\Gs-2<\R s<\Gs\}$, $\Gs\in(0,\Gb)$, $0<\Gb<1$.

{\sl Find two functions, $\hat u_1(\Ga,s)$ and  $\hat u_1(\Ga,s)$, holomorphic everywhere in the strip
$\GP$ except for the point $s=0$ and $s=\nu-1$, where they have  simple poles, vanishing at the infinite points
$\tau\pm i\infty$, $\Gs-2\le\tau\le \Gs$, H\"older-continuous in the strip
up to the boundary $\GO=\{\R s=\Gs\}$ and  $\GO_{-2}=\{\R s=\Gs-2\}$ and satisfying the boundary conditions in the contour $\GO$
 $$
 \hat u_1(\Ga,s)+\fr{i[(a_d^2-a_s^2)(s-1)-\nu a_s^2]}{(1-a_d^2)\Ga}\hat u_2(\Ga,s-1)+\fr{a_s^2(s-2)(s-1-\nu)}{(1-a_d^2)\Ga^2}
 \hat u_1(\Ga,s-2)
$$
$$
 =-\fr{\hat f_1(\Ga,s)}{(1-a_d^2)V^2\Ga^2},  \quad s\in\GO,
$$
$$
 \hat u_2(\Ga,s)+\fr{i[(a_d^2-a_s^2)(s-1)-\nu (a_d^2-2a_s^2)]}{(1-a_s^2)\Ga}\hat u_1(\Ga,s-1)+\fr{a_d^2(s-2)(s-1-\nu)}{(1-a_s^2)\Ga^2}
 \hat u_2(\Ga,s-2)
 $$
 \beq
 =-\fr{\hat f_2(\Ga,s)}{(1-a_s^2)V^2\Ga^2}, \quad s\in\GO.
\label{3.13}
\eeq
}

We further wish to make the coefficients of the functions $\hat u_1(\Ga,s-2)$ and $\hat u_1(\Ga,s-2)$ equal 1 and at the same time not to change the coefficients of the functions  $\hat u_1(\Ga,s)$ and $\hat u_1(\Ga,s)$.  
This is possible to achieve by factorizing the coefficients of the functions $\hat u_j(\Ga,s-2)$ and introducing two new functions
\beq
\GF_j(\Ga, s)=\fr{\GG(1-\fr{s}{2})|\Ga|^s}{\GG(\fr{s+1-\nu}{2})2^s\Gb_j^{s/2}}\hat u_j(\Ga,s), \quad j=1,2,
\label{3.14}
\eeq
where
\beq
\Gb_1=\fr{a_s^2}{a_d^2-1}, \quad \Gb_2=\fr{a_d^2}{a_s^2-1}.
\label{3.15}
\eeq
 In what follows we confine ourselves to the subsonic regime that is assume that $V<c_s$ and therefore both of the parameters
 $\Gb_1$ and $\Gb_2$ are positive.
In the strip $\GP$, these functions have a simple pole at the point $s=0$ and a removable singularity at the point $s=\nu-1$.
The new Carleman boundary value problem is simpler and is stated as follows.

{\sl Find two functions, $\GF_1(\Ga, s)$ and $\GF_2(\Ga, s)$, holomorphic everywhere in the strip
$\GP$ except for the point $s=0$, where they have a simple pole, bounded at the infinite points
$\tau\pm i\infty$, $\Gs-2\le\tau\le \Gs$, H\"older-continuous in the strip
up to the boundary and satisfying the boundary conditions in the contour $\GO$
$$
\GF_1(\Ga,s)-i\sgn\Ga G_1(s)\GF_2(\Ga,s-1)+ \GF_1(\Ga,s-2)=g_1(\Ga,s), 
$$
\beq
\GF_2(\Ga,s)-i\sgn\Ga G_2(s)\GF_1(\Ga,s-1)+ \GF_2(\Ga,s-2)=g_2(\Ga,s),\quad s\in\GO.
\label{3.16}
\eeq
The coefficients of the problem are given by
$$
G_j(s)=b_j(s)
\fr{\GG(1-\fr{s}{2})\GG(\fr{s-\nu}{2})}
{\GG(\fr{s+1-\nu}{2})\GG(\fr{3-s}{2})}, \quad j=1,2,
$$
\beq
b_1(s)=\fr{(a_d^2-a_s^2)(s-1)-\nu a_s^2}{2(a_d^2-1)\Gb_2^{1/2-s/2}\Gb_1^{s/2}},\quad 
\quad
b_2(s)=\fr{(a_d^2-a_s^2)(s-1)-\nu (a_d^2-2a_s^2)}{2(a_s^2-1)\Gb_1^{1/2-s/2}\Gb_2^{s/2}}.
\label{3.17'}
\eeq
The right-hand sides are
\beq
g_1(\Ga,s)=\fr{|\Ga|^{s-2}\GG(1-\fr{s}{2})\hat f_1(\Ga,s)}{V^2(a_d^2-1)2^s\Gb_1^{s/2}\GG(\fr{s+1-\nu}{2})},
\quad
g_2(\Ga,s)=\fr{|\Ga|^{s-2}\GG(1-\fr{s}{2})\hat f_2(\Ga,s)}{V^2(a_s^2-1)2^s\Gb_2^{s/2}\GG(\fr{s+1-\nu}{2})}.
\label{3.17}
\eeq
}

On moving the middle terms $i\sgn\Ga G_1(s)\GF_2(\Ga,s-1)$ and $i\sgn\Ga G_2(s)\GF_1(\Ga,s-1)$
to the right-hand sides and taking into account that the functions $\GF_1(\Ga,s)$ 
and  $\GF_2(\Ga,s)$ have simple pole at the point $s=0\in\GP$  
we  write the general representations of the functions $\GF_j(\Ga,s)$ in the interior of the strip 
\beq
 \GF_j(\Ga, s)=\fr{1}{4i}\int_\GO\fr{g_j(\Ga,p)+i\sgn\Ga G_j(p)\GF_{3-j}(\Ga,p-1)}{\sin\fr{\pi}{2}(p-s)}dp+\fr{C_j(\Ga)}{\sin\fr{\pi s}{2}}, \quad s\in\GP,
 \quad j=1,2.
 \label{3.18}
 \eeq
Here, $C_1(\Ga)$ and $C_2(\Ga)$ are arbitrary functions of $\Ga$.
On  the contours $\GO$ and $\GO_{-2}$,
 by the Sokhotski-Plemelj formulas for a strip,  the functions $\GF_j(\Ga,s)$ 
 have the form
$$
\GF_j(\Ga,s_0)=\fr{1}{2}[g_j(\Ga,s_0)+i\sgn\Ga G_j(s_0)\GF_{3-j}(\Ga,s_0-1)] +\fr{C_j(\Ga)}{\sin\fr{\pi s_0}{2}}+ \CJ_j(\Ga, s_0), \quad s_0\in\GO,
$$
\beq
\GF_j(\Ga,s_0-2)=\fr{1}{2}[g_j(\Ga,s_0)+i\sgn\Ga G_j(s_0)\GF_{3-j}(\Ga,s_0-1)]- \fr{C_j(\Ga)}{\sin\fr{\pi s_0}{2}}- \CJ_j(\Ga, s_0), \quad s_0-2\in\GO_{-2},
\label{3.19}
\eeq
where $\CJ_j(\Ga, s_0)$ is the Cauchy principal value of the integral
\beq
\CJ_j(\Ga,s_0)=\fr{1}{4i}\int_\GO\fr{g_j(\Ga,p)+i\sgn\Ga G_j(p)\GF_{3-j}(\Ga,p-1)}{\sin\fr{\pi}{2}(p-s_0)}dp, \quad s_0\in\GO, \quad j=1,2.
\label{3.20}
\eeq
The functions $\GF_{3-j}(\Ga,p-1)$ and the functions $C_j(\Ga)$  in the representation formula (\ref{3.18}) are unknown.
On replacing  $s$ by $s-1$, $s\in\GO$, we arrive at the following system of two integral equations for the unknown functions:
$$
\GF_1(\Ga,s-1)=\fr{1}{4i}\int_\GO\fr{g_1(\Ga,p)+i\sgn\Ga G_1(p)\GF_2(\Ga,p-1)}{\cos\fr{\pi}{2}(p-s)}dp-\fr{C_1(\Ga)}{\cos\fr{\pi s}{2}},
$$
\beq
\GF_2(\Ga,s-1)=\fr{1}{4i}\int_\GO\fr{g_2(\Ga,p)+i\sgn\Ga G_2(p)\GF_1(\Ga,p-1)}{\cos\fr{\pi}{2}(p-s)}dp-\fr{C_2(\Ga)}{\cos\fr{\pi s}{2}}, \quad s\in\GO.
\label{3.21}
\eeq

\setcounter{equation}{0}

\section{Functions $C_1(\Ga)$ and $C_2(\Ga)$}

The arbitrary functions $C_1(\Ga)$ and $C_2(\Ga)$ are to be fixed by satisfying the boundary conditions (\ref{3.3}).
To do this we express the Fourier-Mellin transforms of the displacements through the functions
$\GF_j(\Ga,s)$
\beq
\hat u_j(\Ga,s)=\fr{\GG(\fr{s+1-\nu}{2})2^s\Gb_j^{s/2}}{\GG(1-\fr{s}{2})|\Ga|^s}\GF_j(\Ga,s), \quad s\in\GP\cup\GO\cup\GO_{-2}.
\label{4.1}
\eeq
On inverting the Mellin transform we have
\beq
\tilde u_j(\Ga,y)=\fr{1}{2\pi i}\int_\GO \fr{\GG(\fr{s+1-\nu}{2})2^s\Gb_j^{s/2}}{\GG(1-\fr{s}{2})|\Ga|^s}\GF_j(\Ga,s) y^{-s}ds.
\label{4.2}
\eeq
If we apply the Cauchy theorem and use the theory of residues, we find
$$
\tilde u_j(\Ga,y)=\fr{1}{2\pi i}\int_{\GO_{-2}} \fr{\GG(\fr{s+1-\nu}{2})2^s\Gb_j^{s/2}}{\GG(1-\fr{s}{2})|\Ga|^s}\GF_j(\Ga,s) y^{-s}ds
$$
\beq
+\left(\mathop{\rm res}\limits_{s=0}+\mathop{\rm res}\limits_{s=\nu-1}\right)\fr{\GG(\fr{s+1-\nu}{2})2^s\Gb_j^{s/2}}{\GG(1-\fr{s}{2})|\Ga|^s}\GF_j(\Ga,s) y^{-s}.
\label{4.3}
\eeq
We compute the residues and derive the following representations of the Fourier transforms of the displacements:
$$
\tilde u_j(\Ga,y)=
\fr{2}{\pi}\GG\left(\fr{1-\nu}{2}\right)C_j(\Ga)
+\fr{2^\nu\Gb_j^{(\nu-1)/2}\GF_j(\Ga,\nu-1)y^{1-\nu}}{|\Ga|^{\nu-1}\GG(\fr{3-\nu}{2})}
$$
\beq
+\fr{1}{2\pi i}\int_{\GO_{-2}} \fr{\GG(\fr{s+1-\nu}{2})2^s\Gb_j^{s/2}}{\GG(1-\fr{s}{2})|\Ga|^s}\GF_j(\Ga,s) y^{-s}ds, \quad y>0.
\label{4.4}
\eeq
Differentiating this expression and evaluating the limit as $y\to 0^+$ give
\beq
\lim_{y\to 0^+} y^\nu\fr{d}{dy}\tilde u_j(\Ga,y)=\fr{2^{\nu+1}\Gb_j^{(\nu-1)/2}\GF_j(\Ga,\nu-1)}{|\Ga|^{\nu-1}\GG(\fr{1-\nu}{2})}.
\label{4.5}
\eeq
Upon substituting these limits into relations (\ref{3.3})
we obtain two equations to be used to fix the functions $C_1(\Ga)$ and $C_2(\Ga)$. They are
\beq
\GF_1(\Ga,\nu-1)=|\Ga|^{\nu-1}\Gg_1\tilde h_1(\Ga), \quad \GF_2(\Ga,\nu-1)=|\Ga|^{\nu-1}\Gg_2\tilde h_2(\Ga),
\label{4.6}
\eeq
where
\beq
\Gg_1=\fr{\GG(\fr{1-\nu}{2})}{\mu_0 2^{\nu+1}\Gb_1^{(\nu-1)/2}}, \quad 
\Gg_2=\fr{\GG(\fr{1-\nu}{2})}{(\Gl_0+2\mu_0) 2^{\nu+1}\Gb_2^{(\nu-1)/2}}.
\label{4.7}
\eeq
To determine the functions $C_1(\Ga)$ and $C_2(\Ga)$, we represent the unknown functions $\GF_1(\Ga,s-1)$ and 
$\GF_2(\Ga,s-1)$ as
\beq
\GF_j(\Ga,s-1)=\GF_j^{(0)}(\Ga,s-1)+C_1(\Ga)\GF_j^{(1)}(\Ga,s-1)+C_2(\Ga)\GF_j^{(2)}(\Ga,s-1), \quad j=1,2,
\label{4.8}
\eeq
and substitute these representations into the system of integral equations (\ref{3.21}). We have
three new systems with the same kernels but different right-hand sides
$$
\GF_j^{(0)}(\Ga,s-1)-\fr{\sgn\Ga}{4}\int_\GO\fr{G_j(p)\GF^{(0)}_{3-j}(\Ga,p-1)dp}{\cos\fr{\pi}{2}(p-s)}
=\fr{1}{4i}\int_\GO\fr{g_j(\Ga,p)dp}{\cos\fr{\pi}{2}(p-s)}, \quad s\in\GO, \quad j=1,2, 
$$
\beq
\GF_j^{(m)}(\Ga,s-1)-\fr{\sgn\Ga}{4}\int_\GO\fr{G_j(p)\GF^{(m)}_{3-j}(\Ga,p-1)dp}{\cos\fr{\pi}{2}(p-s)}
=-\fr{\Gd_{jm}}{\cos\fr{\pi s}{2}}, \quad s\in\GO, \quad j=1,2, \quad m=1,2.
\label{4.9}
\eeq
Note that the new systems do not possess the functions   $C_1(\Ga)$ and $C_2(\Ga)$. From (\ref{4.6}) these functions 
solve the following system of two equations
\beq
\GF_j^{(1)}(\Ga,\nu-1)C_1(\Ga)+\GF_j^{(2)}(\Ga,\nu-1)C_2(\Ga)=|\Ga|^{\nu-1}\Gg_j\tilde h_j(\Ga)-\GF_j^{(0)}(\Ga,\nu-1), \quad j=1,2.
\label{4.10}
\eeq
and have the form
$$
C_1(\Ga)=\fr{1}{\GD(\Ga)}\{[|\Ga|^{\nu-1}\Gg_1\tilde h_1(\Ga)-\GF_1^{(0)}(\Ga,\nu-1)]\GF_2^{(2)}(\Ga,\nu-1)
$$$$
-[|\Ga|^{\nu-1}\Gg_2\tilde h_2(\Ga)-\GF_2^{(0)}(\Ga,\nu-1)]\GF_1^{(2)}(\Ga,\nu-1)\},
$$
$$
C_2(\Ga)=\fr{1}{\GD(\Ga)}\{[|\Ga|^{\nu-1}\Gg_2\tilde h_2(\Ga)-\GF_2^{(0)}(\Ga,\nu-1)]\GF_1^{(1)}(\Ga,\nu-1)
$$
\beq
-[|\Ga|^{\nu-1}\Gg_1\tilde h_1(\Ga)-\GF_1^{(0)}(\Ga,\nu-1)]\GF_2^{(1)}(\Ga,\nu-1)\},
\label{4.11}
\eeq
where
\beq
\GD(\Ga)=\GF_1^{(1)}(\Ga,\nu-1)\GF_2^{(2)}(\Ga,\nu-1)-\GF_1^{(2)}(\Ga,\nu-1)\GF_2^{(1)}(\Ga,\nu-1).
\label{4.12}
\eeq

\setcounter{equation}{0}

\section{Point force running along the boundary}

Assume now that there are no body forces, and a point force $\BF=(H_1,H_2)$ applied at a point $\Gx_0$ is
 running along the boundary of the half-plane at subsonic speed $V$ that is
\beq
f_j(\Gx,y)= 0, \quad y>0, \quad  h_j(\Gx)=H_j\Gd(\Gx-\Gx_0), \quad |\Gx|<\infty, \quad j=1,2.
\label{5.1}
\eeq
Then $\GF_j^{(0)}(s-1)=0$, $\tilde h_j(\Ga)=H_j e^{i\Ga\Gx_0}$, and the functions $C_1(\Ga)$ and $C_2(\Ga)$ are given by
$$
C_1(\Ga)=\fr{|\Ga|^{\nu-1} e^{i\Ga\Gx_0}}{\GD(\Ga)}[
H_1\Gg_1\GF_2^{(2)}(\Ga,\nu-1)-H_2\Gg_2\GF_1^{(2)}(\Ga,\nu-1)],
$$
\beq
C_2(\Ga)=\fr{|\Ga|^{\nu-1} e^{i\Ga\Gx_0}}{\GD(\Ga)}[
H_2\Gg_2\GF_1^{(1)}(\Ga,\nu-1)-H_1\Gg_1\GF_2^{(1)}(\Ga,\nu-1)].
\label{5.2}
\eeq

\subsection{Free of $\Ga$ integral equations}

The two nonzero functions $\GF_j^{(m)}(\Ga,s-1)$ ($j,m=1,2$) depend  on $\Ga$ only because of the presence of $\sgn\Ga$
in the governing system of integral equations (\ref{4.9}).
It will be convenient to introduce new functions associated with $\GF_j^{(m)}(s)$ as follows:
$$
\GF_j(\Ga,s)=C_1(\Ga)\GF_j^{(1)}(\Ga,s)+C_2(\Ga)\GF_j^{(2)}(\Ga,s)
$$
\beq
=|\Ga|^{\nu-1}e^{i\Ga\Gx_0}
\left\{\begin{array}{cc}
C_{1+}\GF_{j+}^{(1)}(s)+C_{2+}\GF_{j+}^{(2)}(s), & \Ga>0,\\
C_{1-}\GF_{j-}^{(1)}(s)+C_{2-}\GF_{j-}^{(2)}(s), & \Ga<0.\\
\end{array}
\right.
\label{5.3}
\eeq
Here, $C_{j\pm}$ are constants and the functions $\GF_{j\pm}^{(m)}(s)$ are independent of $\Ga$. On comparing the relations (\ref{5.2}) and
(\ref{5.3}) we find the constants $C_{1\pm}$ and $C_{2\pm}$
$$
C_{1\pm}=\fr{1}{\GD_\pm}[H_1\Gg_1\GF_{2\pm}^{(2)}(\nu-1)-H_2\Gg_2\GF_{1\pm}^{(2)}(\nu-1)],
$$
\beq
C_{2\pm}=\fr{1}{\GD_\pm}[H_2\Gg_2\GF_{1\pm}^{(1)}(\nu-1)-H_1\Gg_1\GF_{2\pm}^{(1)}(\nu-1)],
\label{5.4}
\eeq
where $\GD_\pm$ are independent of $\Ga$ 
\beq
\GD_\pm=\GF_{1\pm}^{(1)}(\nu-1)\GF_{2\pm}^{(2)}(\nu-1)-\GF_{1\pm}^{(2)}(\nu-1)\GF_{2\pm}^{(1)}(\nu-1).
\label{5.5}
\eeq
As for the functions $\GF_{j\pm}^{(m)}(s)$, $s\in\GP$, they are expressed through the functions $\GF^{(m)}_{3-j\pm}(p-1)$,
$\Gs\in\GO$, as  
\beq
 \GF^{(m)}_{j\pm}(s)=\pm\fr{1}{4}\int_\GO\fr{G_j(p)\GF^{(m)}_{3-j\pm}(p-1)}{\sin\fr{\pi}{2}(p-s)}dp+\fr{\Gd_{jm}}{\sin\fr{\pi s}{2}}, \quad s\in\GP,
 \quad j,m=1,2,
 \label{5.6}
 \eeq
while the functions  $\GF^{(m)}_{j\pm}(s-1)$, $s\in\GO$, solve the system of integral equations 
\beq
\GF_{j\pm}^{(m)}(s-1)=\pm\fr{1}{4}\int_\GO\fr{G_j(p)\GF^{(m)}_{3-j\pm}(p-1)dp}{\cos\fr{\pi}{2}(p-s)}
-\fr{\Gd_{jm}}{\cos\fr{\pi s}{2}}, \quad s\in\GO, \quad j,m=1,2.
\label{5.7}
\eeq

\subsection{Integral representations and asymptotic expansions}

Based on the representation (\ref{5.3}) we aim to derive formulas for the displacements. The Fourier-Mellin double integral transformation of the displacements  (\ref{5.3}) can be written in another, equivalent, form
\beq
\hat u_j(\Ga,s)=\fr{e^{i\Ga\Gx_0}}{|\Ga|^{s+1-\nu}}\left(\fr{\Psi_{j+}(s)+\Psi_{j-}(s)}{2}+
\fr{\Psi_{j+}(s)-\Psi_{j-}(s)}{2}\sgn\Ga\right),
\label{5.8}
\eeq
where
\beq
\Psi_{j\pm}(s)=\fr{\GG(\fr{s+1-\nu}{2})2^s\Gb_j^{s/2}}{\GG(1-\fr{s}{2})}[C_{1\pm}\GF_{j\pm}^{(1)}(s)+ C_{2\pm}\GF_{j\pm}^{(2)}(s)].
\label{5.9}
\eeq
On inverting the Fourier and Mellin transforms we derive integral representations of the displacements. They are
$$
u_j(\Gx,y)=\fr{1}{4\pi^2 i}\left(\int_\GO\fr{\Psi_{j+}(s)+\Psi_{j-}(s)}{2}\CI_0(\Gx_0-\Gx,s)y^{-s}ds
\right.
$$
\beq
\left.
+\int_\GO\fr{\Psi_{j+}(s)-\Psi_{j-}(s)}{2}\CI_1(\Gx_0-\Gx,s)y^{-s}ds\right),
\label{5.10}
\eeq
where
\beq
\CI_0(\Gx)=\int_{-\infty}^\infty\fr{e^{i\Ga \Gx}d\Ga}{|\Ga|^{s+1-\nu}}, \quad 
\CI_1(\Gx)=\int_{-\infty}^\infty\fr{e^{i\Ga \Gx}\sgn\Ga d\Ga}{|\Ga|^{s+1-\nu}}.
\label{5.11}
\eeq 
These two integrals can be explicitly evaluated {\bf(\ref{gra})}, formulas 3.761(9) and  3.761(4). This gives
$$
\CI_0(\Gx)=\fr{2\GG(\nu-s)}{|\Gx-\Gx_0|^{\nu-s}}\cos\fr{\nu-s}{2}\pi,
$$
\beq 
\CI_1(\Gx)=\fr{2i\GG(\nu-s)}{|\Gx-\Gx_0|^{\nu-s}}\sin\fr{\nu-s}{2}\pi\sgn\Gx,
\quad 0<\R s<\nu.
\label{5.12}
\eeq
We substitute our findings into formula (\ref{5.10}) and arrive at the following integral representations of the displacements:
\beq
u_j(\Gx,y)=\fr{1}{4\pi^2 i}\int_\GO\left[\Psi_{j+}(s) e^{i\pi\Gvk(\nu-s)/2}+\Psi_{j-}(s) e^{- i\pi\Gvk(\nu-s)/2}\right]\fr{\GG(\nu-s)y^{-s}ds}{|\Gx-\Gx_0|^{\nu-s}}, \quad y>0, \quad |\Gx|<\infty,
\label{5.13}
\eeq
where $\Gvk=\sgn(\Gx_0-\Gx)$. By applying the Cauchy theorem we shift the contour of integration $\GO$ to $\GO_{-2}$.
Upon computing the residues at the simple poles $s=0$ and $s=\nu-1$ we have 
$$
u_j(\Gx,y)=\fr{\GG(\fr{\nu}{2}) 2^{\nu-1}}{\pi^{3/2}\cos\fr{\pi\nu}{2}|\Gx-\Gx_0|^\nu}[e^{i\pi\Gvk\nu/2}(C_{1+}\Gd_{j1}+C_{2+}\Gd_{j2})+
e^{-i\pi\Gvk\nu/2}(C_{1-}\Gd_{j1}+C_{2-}\Gd_{j2})]
$$
$$
+\fr{i\Gvk 2^{\nu-1} \Gb_j^{(\nu-1)/2}y^{1-\nu}}{\pi\GG(\fr{3-\nu}{2})|\Gx-\Gx_0|}[C_{1+}\GF_{j+}^{(1)}(\nu-1)+
C_{2+}\GF_{j+}^{(2)}(\nu-1)-C_{1-}\GF_{j-}^{(1)}(\nu-1)-
C_{2-}\GF_{j-}^{(2)}(\nu-1)]
$$
\beq
+\fr{1}{4\pi^2 i}\int_{\GO_{-2}}\left[\Psi_{j+}(s) e^{i\pi\Gvk(\nu-s)/2}+\Psi_{j-}(s) e^{-i\pi\Gvk(\nu-s)/2}\right]\fr{\GG(\nu-s)y^{-s}ds}{|\Gx-\Gx_0|^{\nu-s}}.
\label{5.14}
\eeq
It is of interest to obtain more terms in the asymptotic expansions. To do this, we continue the solution through
the contour $\GO_{-2}$. Denote the left and right limits of the functions $\GF_{j\pm}^{(m)}(s)$ on the contour $\GO_{-2}$
 by
 $$
 \GF_{j\pm}^{(m)}(s^+)=\lim_{s\to s_0-2^+} \GF_{j\pm}^{(m)}(s)= \GF_{j\pm}^{(m)}(s_0-2),
 $$
 \beq
 \GF_{j\pm}^{(m)}(s^-)=\lim_{s\to s_0-2^-} \GF_{j\pm}^{(m)}(s)= -\GF_{j\pm}^{(m)}(s_0), \quad s_0\in\GO. 
 \label{5.15}
 \eeq
 From the solution representation formulas (\ref{5.7}) we derive the jumps of the functions $\GF_{j\pm}^{(m)}(s)$   
  when $s$ crosses the contour $\GO_{-2}$
 \beq
\GF_{j\pm}^{(m)}(s^-)= \GF_{j\pm}^{(m)}(s^+)\mp iG_j(s+2)\GF_{3-j\pm}^{(m)}(s+1),\quad s\in\GO_{-2}.
 \label{5.16}
 \eeq 
 The functions $\GF_{j\pm}^{(m)}(s^-)$ defined by this relations admit analytic continuation to the left into the strip
$\GP_{-1}=\{\Gs-3<\R s<\Gs-2\}$. They are 
 holomorphic everywhere  in this strip  except for
 the point $s=-2$ where they have simple poles. The functions $G_j(s+2)$ does not have poles
 in the strip $\GP_{-1}$ (the poles $s=\nu-2$ and $s=\nu-4$ are outside this strip).
 The functions $\GF_{j\pm}^{(m)}(s+1)$ have simple poles at the zeros of $\cos\fr{\pi s}{2}$. However none
 of them falls into the strip $\GP_{-1}$. Therefore the only two poles in the strip $\GP_{-1}$
 of the integrand in the integral (\ref{5.14})
  are the points $s=-2$ and $s=\nu-3$.  We utilize the theory of residues again
 and replace the contour $\GO_{-2}$  by  $\GO_{-3}=\{\R s =\Gs-3\}$. Having computed the residues of the integrand at the points
 $s=-2$ and $s=\nu-3$ we rewrite the resulting integral representation as an asymptotic expansion for small
 $\Gn=y/|\Gx-\Gx_0|$
 \beq
 u_j(\Gx,y)=\fr{1}{|\Gx-\Gx_0|^\nu}
 \left(d_{j0}+d_{j1}\Gn^{1-\nu}+d_{j2}\Gn^2+d_{j3}\Gn^{3-\nu}\right)+R^-_j(\Gx-\Gx_0,y),
 \label{5.17}
 \eeq
 where
 $$
 d_{j0}=\fr{\GG(\fr{\nu}{2}) 2^{\nu-1}}{\pi^{3/2}\cos\fr{\pi\nu}{2}}
[e^{i\pi\Gvk\nu/2}(C_{1+}\Gd_{j1}+C_{2+}\Gd_{j2})+
e^{-i\pi\Gvk\nu/2}(C_{1-}\Gd_{j1}+C_{2-}\Gd_{j2})],
$$
$$
d_{j1}=\fr{i\Gvk2^{\nu-1} \Gb_j^{(\nu-1)/2}}{\pi\GG(\fr{3-\nu}{2})}[C_{1+}\GF_{j+}^{(1)}(\nu-1)+
C_{2+}\GF_{j+}^{(2)}(\nu-1)-C_{1-}\GF_{j-}^{(1)}(\nu-1)-
C_{2-}\GF_{j-}^{(2)}(\nu-1)],
$$
$$
d_{j2}=-\fr{\nu d_{j0}}{2\Gb_j},
$$
$$
d_{j3}=\fr{i\Gvk2^{\nu-2} \Gb_j^{(\nu-3)/2}}{\pi\GG(\fr{5-\nu}{2})}[C_{1+}\GF_{j+}^{(1)}(\nu-3)+
C_{2+}\GF_{j+}^{(2)}(\nu-3)-C_{1-}\GF_{j-}^{(1)}(\nu-3)-
C_{2-}\GF_{j-}^{(2)}(\nu-3)],
$$
$$
R^-_j(\Gx-\Gx_0,y)=\fr{1}{4\pi^2 i}\int_{\GO_{-3}}
\fr{\GG(\nu-s)\GG(\fr{s+1-\nu}{2})2^s\Gb_j^{s/2} y^{-s}}{|\Gx-\Gx_0|^{\nu-s}\GG(1-\fr{s}{2})}
\left\{\left[ C_{1+}\left(\GF_{j+}^{(1)}(s)+i G_j(s+2)\GF_{3-j+}^{(1)}(s+1)\right)
\right.
\right.
$$
$$
\left.
+C_{2+}\left(\GF_{j+}^{(2)}(s)+i G_j(s+2)\GF_{3-j+}^{(2)}(s+1)\right)\right]e^{i\pi\Gvk(\nu-s)/2}
$$
$$
+\left[C_{1-}\left(\GF_{j-}^{(1)}(s)-i G_j(s+2)\GF_{3-j-}^{(1)}(s+1)\right)
\right.
$$
\beq
\left.
\left.
+C_{2-}\left(\GF_{j-}^{(2)}(s)-i G_j(s+2)\GF_{3-j-}^{(2)}(s+1)\right)\right]e^{-i\pi\Gvk(\nu-s)/2}
\right\}ds.
\label{5.18}
\eeq 
Due to the boundary conditions, for $|\Gx-\Gx_0|\ne 0$, we expect that the coefficients $d_{j1}=0$. This
is confirmed by numerical tests. Further, because of the periodicity, $d_{j3}=0$ as well. When simplified, the asymptotic 
expansion (\ref{5.17}) reads
\beq
u_j(\Gx,y)\sim \fr{d_{j0}}{|\Gx-\Gx_0|^\nu}
 \left[1-\fr{\nu}{2\Gb_j}\Gn^2+O(\Gn^4)\right], \quad  \Gn=\fr{y}{|\Gx-\Gx_0|} \; {\rm is \; \, small}.
 \label{5.17'}
 \eeq
Formulas (\ref{5.17}) and (\ref{5.17'}) contain the integral $R^-_j(\Gx-\Gx_0,y)$, which may be further transformed to an integral
over the contour $\GO_{-4}=\{\R s =\Gs-4\}$ by continuing meromorphically the integrand to the strip $\GP_{-2}=\{\Gs-3<s<\Gs-4\}$.
This may bring us extra terms in the asymptotic expansions for the displacements for small $\Gn$.

We now wish to derive an asymptotic expansion of the displacement for large $\Gn$. 
On continuing the solution through the contour $\GO$ into the strip $\GP_1=\{\Gs<\R s<\Gs+1\}$ we employ
the relation between the left and right limits $\GF_{j\pm}^{(m)}(s^-)=\GF_{j\pm}^{(m)}(s_0)$ and 
$\GF_{j\pm}^{(m)}(s^+)=-\GF_{j\pm}^{(m)}(s_0-2)$, respectively, 
on the contour $\GO$, $s_0\in\GO$, $s_0^\pm\in\GO^\pm$,
\beq
\GF_{j\pm}^{(m)}(s^+)=\GF_{j\pm}^{(m)}(s^-)\mp i G_j(s_0)\GF_{3-j\pm}^{(m)}(s_0-1), \quad s_0\in\GO.
\label{5.19}
\eeq
The result of this transition is
$$
u_j(\Gx,y)=\fr{1}{4\pi^2 i}\int_\GO\fr{\GG(\nu-s)\GG(\fr{s+1-\nu}{2})2^s\Gb_j^{s/2}y^{-s}}
{\GG(1-\fr{s}{2})|\Gx-\Gx_0|^{\nu-s}}
$$
$$
\times\left\{\left[ C_{1+}\left(\GF_{j+}^{(1)}(s^+)+i G_j(s)\GF_{3-j+}^{(1)}(s-1)\right)
+C_{2+}\left(\GF_{j+}^{(2)}(s^+)+i G_j(s)\GF_{3-j+}^{(2)}(s-1)\right)\right]\right.
$$
$$
\times e^{i\pi\Gvk(\nu-s)/2}+
\left[ C_{1-}\left(\GF_{j+}^{(1)}(s^+)-i G_j(s)\GF_{3-j-}^{(1)}(s-1)\right)
\right.
$$
\beq
\left.
\left.
+C_{2-}\left(\GF_{j-}^{(2)}(s^+)-i G_j(s)\GF_{3-j-}^{(2)}(s-1)\right)\right]e^{-i\pi\Gvk(\nu-s)/2}\right\}ds.
\label{5.20}
\eeq
The integrand has a simple pole at the point $s=1$ (due to $\GF_{3-j+}^{(1)}(s-1)$) and a pole of order 2 at the point
$s=\nu$ (due to $\GG(s-\nu)$ and $G_j(s)$). To avoid the unnecessary complications associated with the second order pole, 
we derive an asymptotic expansion for the derivative $\Md u_j/\Md\Gx$. The integrand of the integral after the differentiation 
has simple zeros at the points $s=1$ and $s=\nu$ in the strip $\GP_1$. On shifting the contour $\GO$ to the right to replace it by $\GO_1=\{\R s=\Gs+1\}$ we obtain
the following representation of the displacements derivative convenient for large $\Gn=y/|\Gx-\Gx_0|$:
\beq
\fr{\Md u_j}{\Md \Gx}(\Gx,y)=\fr{1}{\pi(\Gx-\Gx_0)y^\nu}[e_{j0}+e_{j1}\Gn^{\nu-1}]+R_j^+(\Gx-\Gx_0,y),
\label{5.21}
\eeq
where
$$
e_{j0}=\fr{i2^\nu b_j(\nu)\Gb_j^{\nu/2}}{\GG(\fr{3-\nu}{2})}
[C_{1+}\GF_{3-j+}^{(1)}(\nu-1)+C_{2+}\GF_{3-j+}^{(2)}(\nu-1)
$$
$$
-C_{1-}\GF_{3-j-}^{(1)}(\nu-1)-C_{2-}\GF_{3-j-}^{(2)}(\nu-1)],
$$
\beq
e_{j1}=-\fr{2\Gvk}{\pi} b_j(1)\sqrt{\Gb_j}\GG(\nu)\GG\left(\fr{1-\nu}{2}\right)
[(C_{1+}\Gd_{j2}+C_{2+}\Gd_{j1})e^{i\pi \Gvk\nu/2}+
(C_{1-}\Gd_{j2}+C_{2-}\Gd_{j1})e^{-i\pi \Gvk\nu/2}].
\label{5.22}
\eeq
The functions $b_j(s)$ are given by (\ref{3.17'}), while the functions $R^-_j(\Gx-\Gx_0,y)$ are the integrals
$$
R^-_j(\Gx-\Gx_0,y)=-\fr{1}{4\pi^2 i(\Gx-\Gx_0)}\int_{\GO_1}
\fr{\GG(\nu-s+1)\GG(\fr{s+1-\nu}{2})2^s\Gb_j^{s/2}y^{-s}}
{\GG(1-\fr{s}{2})|\Gx-\Gx_0|^{\nu-s}}
$$
$$
\times\{[C_{1+}(\GF_{j+}^{(1)}(s)+iG_j(s)\GF_{3-j+}^{(1)}(s-1))
+C_{2+}(\GF_{j+}^{(2)}(s)+iG_j(s)\GF_{3-j+}^{(2)}(s-1))]e^{i\pi\Gvk(\nu-s)/2}
$$
\beq
+[C_{1-}(\GF_{j-}^{(1)}(s)-iG_j(s)\GF_{3-j-}^{(1)}(s-1))
+C_{2-}(\GF_{j-}^{(2)}(s)-iG_j(s)\GF_{3-j-}^{(2)}(s-1))]e^{-i\pi\Gvk(\nu-s)/2}\}ds.
\label{5.23}
\eeq
As in the case of small $\Gn$, it is possible to recover more terms in the representation (\ref{5.21}) by continuing the 
solution further to the right.

We conclude this section by writing down asymptotic expansions of the stresses $\Gs_{12}$ and $\Gs_{22}$ for small $\Gn$.
On substituting the representations (\ref{5.17'}) into the formulas
\beq
\Gs_{12}=\mu_0y^\nu\left(\fr{\Md u_1}{\Md y}+\fr{\Md u_2}{\Md \Gx}\right),\quad 
\Gs_{22}=y^\nu\left[\Gl_0\fr{\Md u_1}{\Md \Gx}+(\Gl_0+2\mu_0)\fr{\Md u_2}{\Md y}\right]
\label{5.24}
\eeq
we arrive at the following asymptotic expansions for small  $\Gn=y/|\Gx-\Gx_0|$:
$$
\fr{\Gs_{12}(\Gx,y)}{\mu_0}\sim \fr{\Gn^\nu}{|\Gx-\Gx_0|}[-\nu d_{20}+2d_{12}\Gn-(\nu+2)d_{22}\Gn^2+O(\Gn^3)],
$$
\beq
\fr{\Gs_{22}(\Gx,y)}{\mu_0}\sim \fr{\Gn^\nu}{|\Gx-\Gx_0|}\left[-\fr{\Gl_0 \nu}{\mu_0} d_{10}+\fr{\Gl_0+2\mu_0}{\mu_0}2d_{22}\Gn-\fr{\Gl_0(\nu+2)}{\mu_0}d_{12}\Gn^2+O(\Gn^3)\right].
\label{5.25}
\eeq

\setcounter{equation}{0}

\section{System of integral equations}

In this section we analyze the system of integral equations (\ref{5.7}) and develop a numerical procedure for its solution.

\subsection{Reduction to a system on the interval $(0,1)$ and its analysis}

First we make the substitutions $p=\Gs+i\tau$ and $s=\Gs+it$  and transform the system (\ref{5.7}) to
the form
$$
\GF_{j\pm}^{(m)}(\Gs-1+it)=\pm\fr{i}{4}\int_{-\infty}^\infty\fr{G_j(\Gs+i\tau)\GF_{3-j\pm}^{(m)}(\Gs-1+i\tau)d\tau}{\cosh\fr{\pi}{2}(\tau-t)}-
\fr{\Gd_{jm}}{\cos\fr{\pi}{2}(\Gs+it)}, 
$$
\beq
-\infty<t<\infty, \quad j,m=1,2.
\label{6.1}
\eeq
Consider first the case of positive $t$ and make the substitutions $x=e^{-\pi t}$. 
Split the integral into integrals over the intervals $(-\infty,0)$ and $(0,\infty)$
and  put $y=e^{\pi \tau}$ when $\tau<0$ and $y=e^{-\pi \tau}$ if $\tau>0$. 
Denote 
\beq
\tau^\pm=\Gs-1\pm\fr{i\ln y}{\pi},\quad t^\pm=\Gs-1\pm\fr{i\ln x}{\pi}
\label{6.2}
\eeq
and rename the unknown functions 
$$
\Gvf_{1\pm}^{(m+)}(y)=y^{-1/2}G_2(\tau^+)\GF_{1\pm}^{(m)}(\tau^+),
\quad
\Gvf_{1\pm}^{(m-)}(y)=y^{-1/2}G_2(\tau^-)\GF_{1\pm}^{(m)}(\tau^-),
$$
\beq
\Gvf_{2\pm}^{(m+)}(y)=y^{-1/2}G_1(\tau^+)\GF_{2\pm}^{(m)}(\tau^+),
\quad
\Gvf_{2\pm}^{(m-)}(y)=y^{-1/2}G_1(\tau^-)\GF_{2\pm}^{(m)}(\tau^-).
\label{6.3}
\eeq
Simple and obvious transformations allow us to rewrite the system (\ref{6.1}) of two equations in an infinite interval as a new system
of four equations in the finite interval
$(0,1)$
$$
\pm\fr{2\pi i\Gvf_{1\pm}^{(m-)}(x)}{G_2(\Gs-\fr{i}{\pi}\ln x)}+\int_0^1\fr{\Gvf_{2\pm}^{(m-)}(y)dy}{y+x}+
\int_0^1\fr{\Gvf_{2\pm}^{(m+)}(y)dy}{1+yx}=\mp \Gd_{m1}f(x),
$$
$$
\pm\fr{2\pi i\Gvf_{1\pm}^{(m+)}(x)}{G_2(\Gs+\fr{i}{\pi}\ln x)}+\int_0^1\fr{\Gvf_{2\pm}^{(m+)}(y)dy}{y+x}+
\int_0^1\fr{\Gvf_{2\pm}^{(m-)}(y)dy}{1+yx}=\pm \Gd_{m1}\ov{f(x)},
$$
$$
\pm\fr{2\pi i\Gvf_{2\pm}^{(m-)}(x)}{G_1(\Gs-\fr{i}{\pi}\ln x)}+\int_0^1\fr{\Gvf_{1\pm}^{(m-)}(y)dy}{y+x}+
\int_0^1\fr{\Gvf_{1\pm}^{(m+)}(y)dy}{1+yx}=\mp \Gd_{m2}f(x),
$$
\beq
\pm\fr{2\pi i\Gvf_{2\pm}^{(m+)}(x)}{G_1(\Gs+\fr{i}{\pi}\ln x)}+\int_0^1\fr{\Gvf_{1\pm}^{(m+)}(y)dy}{y+x}+
\int_0^1\fr{\Gvf_{1\pm}^{(m-)}(y)dy}{1+yx}=\pm  \Gd_{m2}\ov{f(x)}, \quad 0<x<1,
\label{6.4}
\eeq
where $m=1,2$ and
\beq
f(x)=\fr{4\pi i}{x e^{\pi i\Gs/2}+ e^{-\pi i\Gs/2}}.
\label{6.5}
\eeq
In fact, equations (\ref{6.4})
represent two independent systems, the ``+" system and ``-" system, of four equations for two sets of four functions, $\Gvf_{1\pm}^{(m-)}(x)$,
$\Gvf_{1\pm}^{(m+)}(x)$, $\Gvf_{2\pm}^{(m-)}(x)$, and $\Gvf_{2\pm}^{(m+)}(x)$. Also, each system needs to be solved twice  for
 two different right-hand sides associated with $m=1$ and $m=2$. 

A numerical algorithm to be developed for the solution of system (\ref{6.4}) has to address two features of the system.
Firstly, the kernel of the first integral in each equation of the system has a fixed singularity at $y=x=0$. Secondly, the 
functions $G_j(\Gs\pm \fr{i}{\pi}\ln x)$ reveal oscillating behavior in a neighborhood of the point $x=0$. 
To understand the nature of this behavior, we analyze
\beq
\fr{\GG(1-\fr{s}{2})\GG(\fr{s-\nu}{2})}
{\GG(\fr{s+1-\nu}{2})\GG(\fr{3-s}{2})}=-\fr{\GG(\fr{s-1}{2})\GG(\fr{s-\nu}{2})}
{\GG(\fr{s+1-\nu}{2})\GG(\fr{s}{2})}\cot\fr{\pi s}{2}
\sim\pm \fr{2i}{s}, \quad s=\Gs+i\tau, \quad \tau\to\pm\infty.
\label{6.6}
\eeq
Consequently formulas (\ref{3.17'}) yield
\beq
G_1(\Gs+i\tau)\sim \pm i\Gl_1\Gb^{(\Gs+i\tau)/2}, \quad G_2(\Gs+i\tau)\sim \pm i\Gl_2\Gb^{-(\Gs+i\tau)/2},\quad  \tau\to\pm\infty,
\label{6.7}
\eeq
where
\beq
\Gb=\fr{\Gb_2}{\Gb_1}, \quad \Gl_1=\fr{a_d^2-a_s^2}{(a_d^2-1)\sqrt{\Gb_2}}, \quad 
  \Gl_2=\fr{a_d^2-a_s^2}{(a_s^2-1)\sqrt{\Gb_1}}.
  \label{6.8}
  \eeq
  This brings us to the relations which describe the oscillatory behavior at the point $x=0$
  of the functions  $G_j(\Gs\pm \fr{i}{\pi}\ln x)$ 
  $$
  G_1(\Gs-\fr{i}{\pi}\ln x)\sim i\Gl_1\Gb^{\Gs/2}x^{-i\Gve}, \quad 
   G_1(\Gs+\fr{i}{\pi}\ln x)\sim -i\Gl_1\Gb^{\Gs/2}x^{i\Gve}, 
   $$
\beq
  G_2(\Gs-\fr{i}{\pi}\ln x)\sim i\Gl_2\Gb^{-\Gs/2}x^{i\Gve}, \quad 
   G_2(\Gs+\fr{i}{\pi}\ln x)\sim -i\Gl_2\Gb^{-\Gs/2}x^{-i\Gve}, \quad x\to 0^+, 
\label{6.9}
\eeq  
  where $\Gve=\fr{1}{2\pi}\ln\Gb$ is a real number. Due to this behavior the solution of the system (\ref{6.4}) also oscillates
  near the point $x=0$, 
  \beq
\Gvf_{j\pm}^{(m-)}(x) \sim  A_{j\pm}^{(m-)} x^{i\Gd_{j}^-},\quad
\Gvf_{j\pm}^{(m+)}(x) \sim A_{j\pm}^{(m+)} x^{i\Gd_{j}^+},\
 \quad x\to 0.
  \label{6.10}
  \eeq
Here, $A_{j\pm}^{(m-)}$ and $A_{j\pm}^{(m+)}$ are complex constants, while  $\Gd_{j}^-$, and  $\Gd_{j}^+$
are real. To determine the parameters $\Gd_{j}^\pm$, we evaluate the singular integral
\beq
\CM(x,\Gd)=\int_0^1\fr{y^{i\Gd}dy}{y+x}, \quad 0<x<1.
\label{6.11}
\eeq
We extend the integral to the interval $0<x<\infty$ and write it as a Mellin-convolution integral
\beq
\CM(x,\Gd)=\int_0^\infty\fr{y^{i\Gd}_-}{1+\fr{x}{y}}\fr{dy}{y}, \quad 0<x<\infty,
  \label{6.12}
  \eeq
  where $y_-^{i\Gd}=0$ if $y>1$ and $y_-^{i\Gd}=y^{i\Gd}$ when $0<y<1$.
  By applying the convolution theorem of the Mellin transform we have
  \beq
  \CM(x,\Gd)=\fr{1}{2 i}\int_{\Gk-i\infty}^{\Gk+i\infty} \fr{x^{-s}ds}{(i\Gd+s)\sin\pi s}, \quad \Gk\in(0,1).
  \label{6.13}
  \eeq 
  This integral can be computed by the theory of residues
  \beq
  \CM(x,\Gd)=\fr{\pi i x^{i\Gd}}{\sinh\pi \Gd}+\sum_{j=0}^\infty\fr{(-1)^jx^j}{i\Gd-j}, \quad 0<x<1.
  \label{6.14}
  \eeq
  Employing this result  and also formulas (\ref{6.10}) we find the behavior of the singular integrals
  in the system (\ref{6.4}). Keeping the oscillating terms and dropping out those which do not oscillate at the point $x=0$
  we obtain from the first and third equations of the system (\ref{6.4})
 $$
 A_{1\pm}^{(m-)}\fr{\Gb^{\Gs/2}x^{i(\Gd_1^- -\Gve)}}{i\Gl_2}\pm A_{2\pm}^{(m-)}\fr{x^{i\Gd_2^-}}{2\sinh \pi\Gd_2^-}\sim \const,
 \quad x\to 0,
 $$
 \beq
 A_{2\pm}^{(m-)}\fr{\Gb^{-\Gs/2}x^{i(\Gd_2^-+\Gve)}}{i\Gl_1}\pm A_{1\pm}^{(m-)}\fr{x^{i\Gd_1^-}}{2\sinh \pi\Gd_1^-}\sim \const, 
 \quad x\to 0.
\label{6.15}
\eeq
 The oscillating terms cancel each other if and only if
 \beq
 \Gd_1^- -\Gd_2^-=\Gve  
  \label{6.16}
  \eeq
  and the determinant of the system 
  $$
  \fr{\Gb^{\Gs/2}}{\Gl_2}A_{1\pm}^{(m-)}\pm\fr{i}{2\sinh\pi\Gd_2^-} A_{2\pm}^{(m-)}=0,
  $$
  \beq
 \fr{\Gb^{-\Gs/2}}{\Gl_1}A_{2\pm}^{(m-)}\pm\fr{i}{2\sinh\pi\Gd_1^-} A_{1\pm}^{(m-)}=0
 \label{6.17}
 \eeq
  with respect to the coefficients $A_{1\pm}^{(m-)}$ and $A_{2\pm}^{(m-)}$
  equals 0 or, equivalently,
  \beq
  \fr{\Gl_1\Gl_2}{4\sinh\pi\Gd_1^-\sinh\pi\Gd_2^-}+1=0, \quad \Gd_2^-= \Gd_1^- -\Gve  .
  \label{6.18}
  \eeq
  If the parameters  $\Gd_1^-$ and $\Gd_2^-$ satisfy equations (\ref{6.18}), then
  the oscillatory terms are cancelled and the coefficients  $A_{1\pm}^{(m-)}$ and $A_{2\pm}^{(m-)}$ are connected by the relation
 \beq
 A_{2\pm}^{(m-)}=\mp\fr{i\Gb^{\Gs/2}\Gl_1}{2\sinh\pi\Gd_1^-}A_{1\pm}^{(m-)}.
 \label{6.19}
 \eeq
There are two sets of solutions of the system of equations  (\ref{6.18}). They are
\beq
\Gd_1^-=\fr{\Gve}{2}\pm l, \quad \Gd_2^-=-\fr{\Gve}{2}\pm l,
\label{6.20}
\eeq
where
\beq
l=\fr{1}{2\pi}\ln(r+\sqrt{r^2-1}), \quad r=\cosh\pi\Gve-\fr{\Gl_1\Gl_2}{2}.
\label{6.21}
\eeq
The oscillatory terms in the second and fourth equations of the system (\ref{6.4})  are analyzed in a similar manner.
The necessary and sufficient conditions for their cancellation read
\beq
\fr{\Gl_1\Gl_2}{4\sinh\pi\Gd_1^+\sinh\pi\Gd_2^+}+1=0, \quad 
\Gd_2^+=\Gd_1^++\Gve.
  \label{6.22}
  \eeq
The analog of the relation (\ref{6.19}) is 
 \beq
 A_{2\pm}^{(m+)}=\pm\fr{i\Gb^{\Gs/2}\Gl_1}{2\sinh\pi\Gd_1^+}A_{1\pm}^{(m+)},
 \label{6.23}
 \eeq
and the sets of solutions of the system (\ref{6.21})  have the form
\beq
\Gd_1^+=-\fr{\Gve}{2}\pm l, \quad \Gd_2^-=\fr{\Gve}{2}\pm l.
\label{6.24}
\eeq
 We choose the following values for the parameters $\Gd_1^\pm$ and $\Gd_2^\pm$:
 \beq
 \Gd_1^-=\Gd_2^+=\fr{\Gve}{2}+l, \quad \Gd_1^+=\Gd_2^-=-\fr{\Gve}{2}+l.
 \label{6.25}
 \eeq

 \subsection{Numerical solution of the system (\ref{6.4})}

For the solution of the system of integral equations (\ref{6.4}) we develop a method based on quadrature formulas 
that counts for the oscillating singularity of the solution at the point $x=0$.
It will be convenient to represent the unknown functions as
\beq
\Gvf_{j\pm}^{(m-)}(x)=x^{i\Gd_j^-}\hat\Gvf_{j\pm}^{(m-)}(x), \quad \Gvf_{j\pm}^{(m+)}(x)=x^{i\Gd_j^+}\hat\Gvf_{j\pm}^{(m+)}(x).
\label{6.26}
\eeq
 We split the interval $[0,1]$ into $N$ subintervals $[x_{k-1},x_k]$ ($k=1,2,\ldots,N$) of the same length $1/N$, $x_k=k/N$, $k=0,1,\ldots,N$, and approximate the unknown functions
 as follows:
 \beq
\hat\Gvf_{j\pm}^{(m-)}(x)=F_{jk\pm}^{(m-)}, \quad \hat\Gvf_{j\pm}^{(m+)}(x)=F_{jk\pm}^{(m+)}, \quad x\in(x_{k-1},x_k], \quad k=1,2,\ldots,N.
\label{6.27}
\eeq
 To approximate the singular integrals in the system (\ref{6.4}), we remove the singularity at the point $x=0$ by writing
 \beq
 \int_0^1\fr{\Gvf_{j\pm}^{(m-)}(y)dy}{y+x_k}=
 \sum_{n=1}^N\int_{x_{n-1}}^{x_n}[\hat\Gf_{j\pm}^{(m-)}(y)-\hat\Gvf_{j\pm}^{(m-)}(0)]
 \fr{y^{i\Gd_j^-}dy}{y+x_k}+\hat\Gvf_{j\pm}^{(m-)}(0)
 \sum_{n=1}^N\int_{x_{n-1}}^{x_n}\fr{y^{i\Gd_j^-}dy}{y+x_k}.
 \label{6.28}
 \eeq
The integrals in the first sum in (\ref{6.28}) are evaluated approximately
\beq
\int_{x_{n-1}}^{x_n}[\hat\Gf_{j\pm}^{(m-)}(y)-\hat\Gvf_{j\pm}^{(m-)}(0)]
 \fr{y^{i\Gd_j^-}dy}{y+x_k}=\fr{F_{jn\pm}^{(m-)}-F_{j1\pm}^{(m-)}}{x_{n-1}+x_k}\fr{x_n^{i\Gd_j^-+1}- x_{n-1}^{i\Gd_j^-+1}}{i\Gd_j^-+1}.
 \label{6.29}
 \eeq
 We  pass to the limit $N\to\infty$ in the second sum and find
 \beq
 \lim_{N\to\infty} \sum_{n=1}^N\int_{x_{n-1}}^{x_n}\fr{y^{i\Gd_j^-}dy}{y+x_k}=
 \int_0^1\fr{y^{i\Gd_j^-} dy}{y+x_k}=\CM(x_k,\Gd_j^-),
 \label{6.30}
 \eeq
 where $\CM(x,\Gd)$ is given by formula (\ref{6.14}).
 On combining our findings we deduce the quadrature formula
 \beq
  \int_0^1\fr{\Gvf_{j\pm}^{(m-)}(y)dy}{y+x_k}=F_{j1\pm}^{(m-)}\CM(x_k,\Gd_j^-)+\sum_{n=2}^N
 \fr{F_{jn\pm}^{(m-)}-F_{j1\pm}^{(m-)}}{x_{n-1}+x_k}\fr{x_n^{i\Gd_j^-+1}- x_{n-1}^{i\Gd_j^-+1}}{i\Gd_j^-+1}.
 \label{6.31}
 \eeq 
 The other integrals in (\ref{6.4}) are regular and their approximation has the form
  \beq
  \int_0^1\fr{\Gvf_{j\pm}^{(m-)}(y)dy}{1+yx_k}=\sum_{n=1}^N
 \fr{F_{jn\pm}^{(m-)}}{1+x_{n-1}x_k}\fr{x_n^{i\Gd_j^-+1}- x_{n-1}^{i\Gd_j^-+1}}{i\Gd_j^-+1}.
 \label{6.32}
 \eeq
 The formulas for the integrals possessing  $\Gvf_{j\pm}^{(m+)}(y)$ coincide with (\ref{6.31}) and (\ref{6.32})  provided
 the upper subscripts with $``-"$ are replaced by the ones with $``+"$.
 By employing these formulas we may approximate the system of singular integral equations (\ref{6.4}) by 
 a linear algebraic system
 $$
\pm\fr{2\pi i F_{1k\pm}^{(m-)}x_k^{i\Gd_1^-}}{G_2(\Gs-\fr{i}{\pi}\ln x_k)}+F_{21\pm}^{(m-)}\CM(x_k,\Gd_1^+)+\sum_{n=2}^N
 \fr{F_{2n\pm}^{(m-)}-F_{21\pm}^{(m-)}}{x_{n-1}+x_k}\fr{x_n^{i\Gd_1^+ +1}- x_{n-1}^{i\Gd_1^+ +1}}{i\Gd_1^+ +1}
 $$
 $$
+\sum_{n=1}^N
 \fr{F_{2n\pm}^{(m+)}}{1+x_{n-1}x_k}\fr{x_n^{i\Gd_1^- +1}- x_{n-1}^{i\Gd_1^- +1}}{i\Gd_1^- +1}
=\mp\Gd_{m1}f(x_k),
$$
 $$
\pm\fr{2\pi i F_{1k\pm}^{(m+)}x_k^{i\Gd_1^+}}{G_2(\Gs+\fr{i}{\pi}\ln x_k)}+F_{21\pm}^{(m+)}\CM(x_k,\Gd_1^-)+\sum_{n=2}^N
 \fr{F_{2n\pm}^{(m+)}-F_{21\pm}^{(m+)}}{x_{n-1}+x_k}\fr{x_n^{i\Gd_1^- +1}- x_{n-1}^{i\Gd_1^- +1}}{i\Gd_1^- +1}
 $$
 $$
+\sum_{n=1}^N
 \fr{F_{2n\pm}^{(m-)}}{1+x_{n-1}x_k}\fr{x_n^{i\Gd_1^+ +1}- x_{n-1}^{i\Gd_1^+ +1}}{i\Gd_1^+ +1}
=\pm\Gd_{m1}\ov{f(x_k)},
$$
 $$
\pm\fr{2\pi i F_{2k\pm}^{(m-)}x_k^{i\Gd_1^+}}{G_1(\Gs-\fr{i}{\pi}\ln x_k)}+F_{11\pm}^{(m-)}\CM(x_k,\Gd_1^-)+\sum_{n=2}^N
 \fr{F_{1n\pm}^{(m-)}-F_{11\pm}^{(m-)}}{x_{n-1}+x_k}\fr{x_n^{i\Gd_1^- +1}- x_{n-1}^{i\Gd_1^- +1}}{i\Gd_1^- +1}
 $$
 $$
+\sum_{n=1}^N
 \fr{F_{1n\pm}^{(m+)}}{1+x_{n-1}x_k}\fr{x_n^{i\Gd_1^+ +1}- x_{n-1}^{i\Gd_1^+ +1}}{i\Gd_1^+ +1}
=\mp\Gd_{m2}f(x_k),
$$ 
  $$
\pm\fr{2\pi i F_{2k\pm}^{(m+)}x_k^{i\Gd_1^-}}{G_1(\Gs+\fr{i}{\pi}\ln x_k)}+F_{11\pm}^{(m+)}\CM(x_k,\Gd_1^+)+\sum_{n=2}^N
 \fr{F_{1n\pm}^{(m+)}-F_{11\pm}^{(m+)}}{x_{n-1}+x_k}\fr{x_n^{i\Gd_1^+ +1}- x_{n-1}^{i\Gd_1^+ +1}}{i\Gd_1^+ +1}
 $$
 \beq
+\sum_{n=1}^N
 \fr{F_{1n\pm}^{(m-)}}{1+x_{n-1}x_k}\fr{x_n^{i\Gd_1^- +1}- x_{n-1}^{i\Gd_1^- +1}}{i\Gd_1^- +1}
=\pm\Gd_{m2}\ov{f(x_k)},\quad k=1,2,\ldots,N.
\label{6.33}
\eeq
For numerical purposes it is helpful to clarify the structure of the matrix of the algebraic system.
Denote it by $A=\{a_{kn}\}$, $k,n=1,2,\ldots 4N$. It may be split into 16 blocks of dimension $N\times N$
\beq
A=\left(
\begin{array}{cccc}
A_{11} & 0 & A_{13} & A_{14}\\
0 & A_{22} & A_{23} & A_{24}\\
A_{24} & A_{23} & A_{33} & 0\\
A_{14} & A_{13} & 0 & A_{44}\\
\end{array}
\right)
\label{6.34}
\eeq
with four blocks being zero matrices. The diagonal blocks $A_{kk}$ are diagonal matrices with the diagonal entries 
$$
a_{kk}=\pm\fr{2\pi i x_k^{i\Gd_1^-}}{G_2(\Gs-\fr{i}{\pi}\ln x_k)}, \quad a_{k+N\,k+N}=\pm\fr{2\pi i x_k^{i\Gd_1^+}}{G_2(\Gs+\fr{i}{\pi}\ln x_k)},
$$
\beq
a_{k+2N\,k+2N}=\pm\fr{2\pi i x_k^{i\Gd_1^+}}{G_1(\Gs-\fr{i}{\pi}\ln x_k)}, \quad 
a_{k+3N\,k+3N}=\pm\fr{2\pi i x_k^{i\Gd_1^-}}{G_1(\Gs+\fr{i}{\pi}\ln x_k)},
\quad k=1,2,\ldots,N.
\label{6.35}
\eeq
Out of the eight blocks left only four are distinct, $A_{13}=A_{42}$, $A_{14}=A_{41}$, $A_{23}=A_{32}$, and $A_{24}=A_{31}$.
The entries of the blocks $A_{13}$ and $A_{24}$  are associated with the singular integrals in the system (\ref{6.4}). They are 
$$
a_{k\,2N+1}=\CM(x_k,\Gd_1^+)-\sum_{n=2}^N 
 \fr{x_n^{i\Gd_1^+ +1}- x_{n-1}^{i\Gd_1^+ +1}}{(i\Gd_1^+ +1)(x_{n-1}+x_k)},
$$
$$
a_{k+N\,3N+1}=\CM(x_k,\Gd_1^-)-\sum_{n=2}^N
 \fr{x_n^{i\Gd_1^- +1}- x_{n-1}^{i\Gd_1^- +1}}{(i\Gd_1^- +1)(x_{n-1}+x_k)},
$$
\beq
a_{k\, 2N+n}= \fr{x_n^{i\Gd_1^+ +1}- x_{n-1}^{i\Gd_1^+ +1}}{(i\Gd_1^+ +1)(x_{n-1}+x_k)},\quad 
a_{k+N\, 3N+n}= \fr{x_n^{i\Gd_1^- +1}- x_{n-1}^{i\Gd_1^- +1}}{(i\Gd_1^- +1)(x_{n-1}+x_k)},\quad n=2,3,\ldots,N.
\label{6.36}
\eeq
 The regular integrals in the system (\ref{6.4}) generate the entries of the blocks $A_{14}$ and $A_{23}$
\beq
a_{k\, 3N+n}=\fr{x_n^{i\Gd_1^- +1}- x_{n-1}^{i\Gd_1^- +1}}{(i\Gd_1^- +1)(1+x_{n-1}x_k)} 
 \quad 
 a_{N+k\,2 N+n}=\fr{x_n^{i\Gd_1^+ +1}- x_{n-1}^{i\Gd_1^+ +1}}{(i\Gd_1^+ +1)(1+x_{n-1}x_k)}, \quad n=1,2,\ldots,N.
 \label{6.37}
\eeq
On introducing new vectors of unknowns by the relations 
\beq
\CF_k=F_{1k\pm}^{(m-)}, \quad \CF_{N+k}=F_{1k\pm}^{(m+)},\quad 
\CF_{2N+k}=F_{2k\pm}^{(m-)},\quad \CF_{3N+k}=F_{2k\pm}^{(m+)},
\label{6.38}
\eeq
and  new vectors for the right-hand side
$$
r_k=\mp\Gd_{m1}f(x_k), \quad r_{N+k}=\pm\Gd_{m1}\ov{f(x_k)},
$$
\beq
r_{2N+k}=\mp\Gd_{m2}f(x_k), \quad r_{2N+k}=\pm\Gd_{m2}\ov{f(x_k)}, \quad k=1,2,\ldots,N,
\label{6.39}
\eeq
we can write the algebraic system in the form
\beq
\sum_{n=1}^{4N} a_{kn}\CF_n=r_k, \quad k=1,2,\ldots,4N.
\label{6.40}
\eeq

\setcounter{equation}{0}

\section{Numerical results}

\begin{figure}[t]
\centerline{
\scalebox{0.6}{\includegraphics{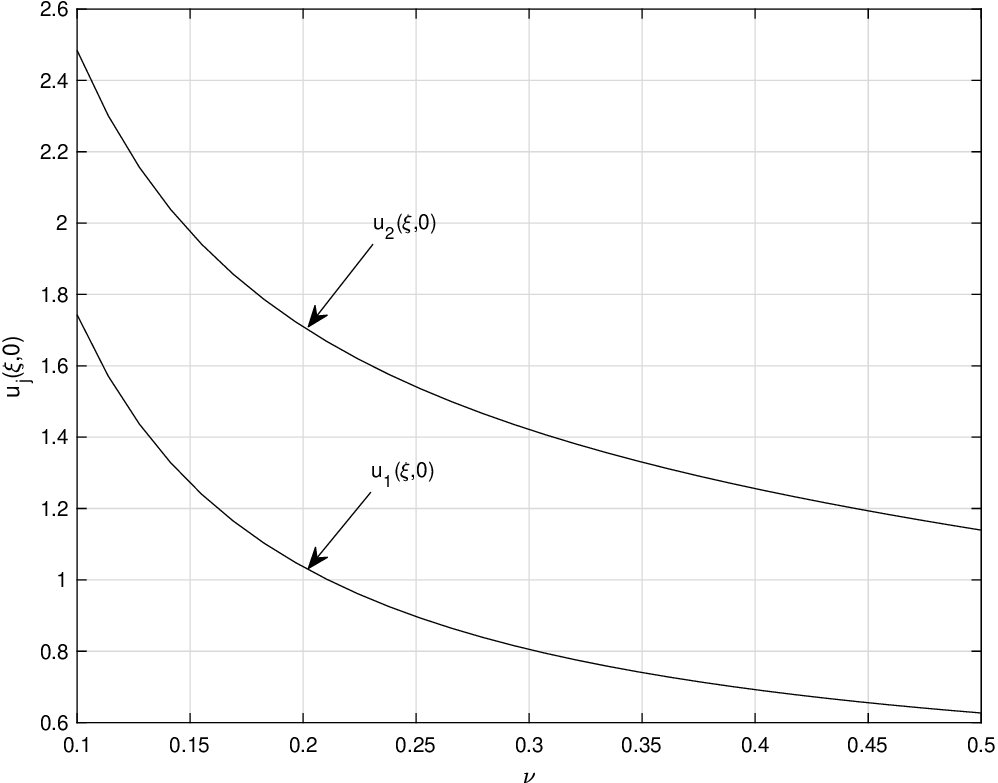}}}
\caption{
The variation of the  tangential and normal displacements  $u_j(\Gx,0)$ on the boundary  $y=0$ when  $0.1\le \nu\le 0.5$,   $\mu_0^{-1}H_1=\mu_0^{-1}H_2=-1$,
$\Gx_0=0$, $\Gx=-1$, $\nu_P=0.3$, $V/c_s=0.2$.}
\label{fig2}
\end{figure}

\begin{figure}[t]
\centerline{
\scalebox{0.6}{\includegraphics{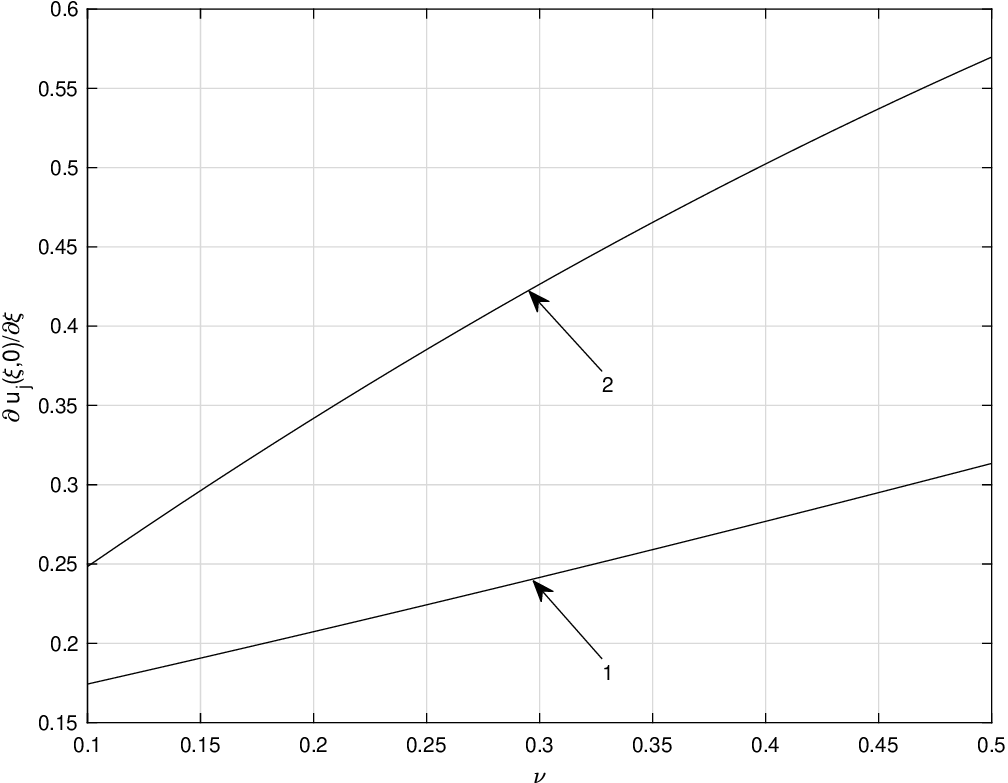}}}
\caption{
The variation of the derivatives of the tangential (curve 1)  and normal (curve 2) displacements  $\fr{\Md}{\Md\xi}u_j(\Gx,0)$ on the boundary  $y=0$ when  $0.1\le \nu\le 0.5$,   $\mu_0^{-1}H_1=\mu_0^{-1}H_2=-1$, 
$\Gx_0=0$, $\Gx=-1$, $\nu_P=0.3$, $V/c_s=0.2$.}
\label{fig3}
\end{figure}

\begin{figure}[t]
\centerline{
\scalebox{0.6}{\includegraphics{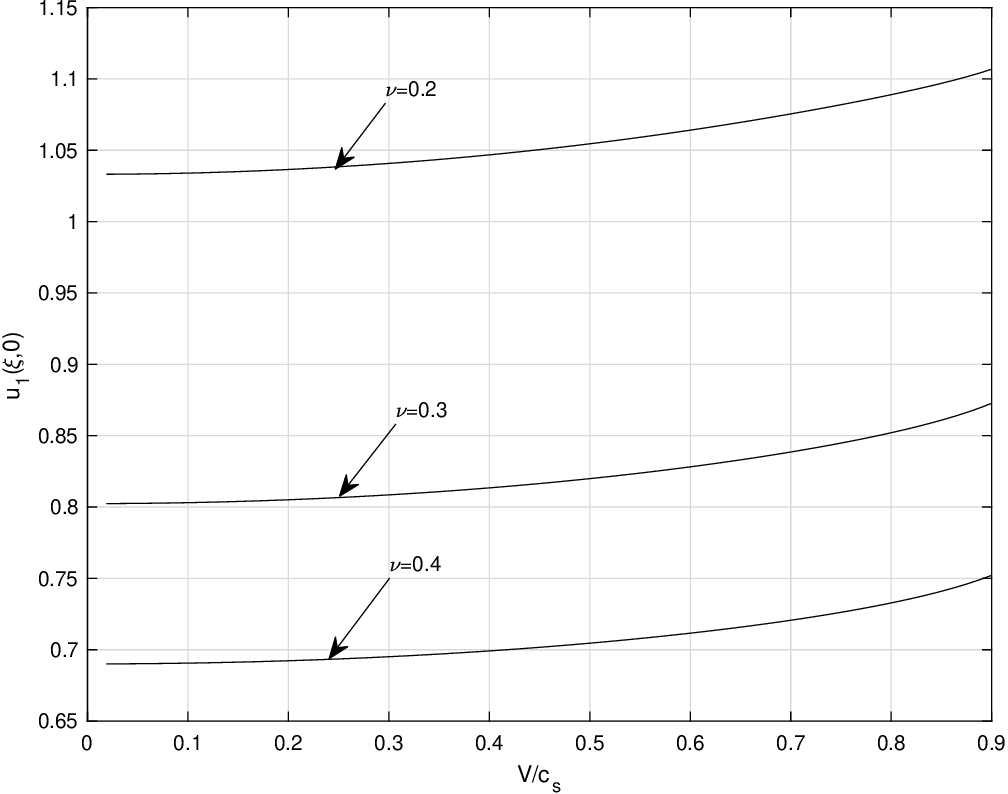}}}
\caption{
The tangential displacements  $u_1(\Gx,0)$ on the boundary  $y=0$ versus  the dimensionless  speed $V/c_s$  
the values $0.2, 0.3, 0.4$ of the parameter $\nu$
when  $\mu_0^{-1}H_1=\mu_0^{-1}H_2=-1$, 
$\Gx_0=0$, $\Gx=-1$, $\nu_P=0.3$.}
\label{fig4}
\end{figure}

\begin{figure}[t]
\centerline{
\scalebox{0.6}{\includegraphics{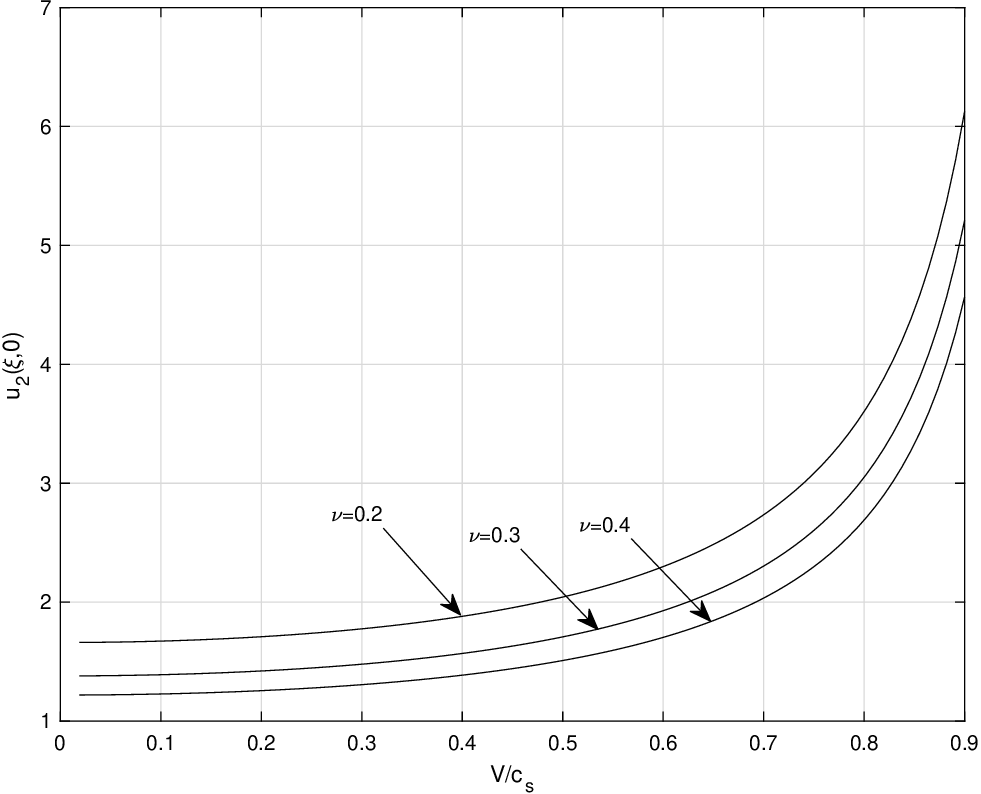}}}
\caption{
The normal displacements  $u_2(\Gx,0)$ on the boundary  $y=0$ 
versus  the dimensionless speed $V/c_s$  for the values $0.2, 0.3, 0.4$ of the parameter $\nu$
when   $\mu_0^{-1}H_1=\mu_0^{-1}H_2=-1$, 
$\Gx_0=0$, $\Gx=-1$, $\nu_P=0.3$.}
\label{fig5}
\end{figure}

\begin{figure}[t]
\centerline{
\scalebox{0.6}{\includegraphics{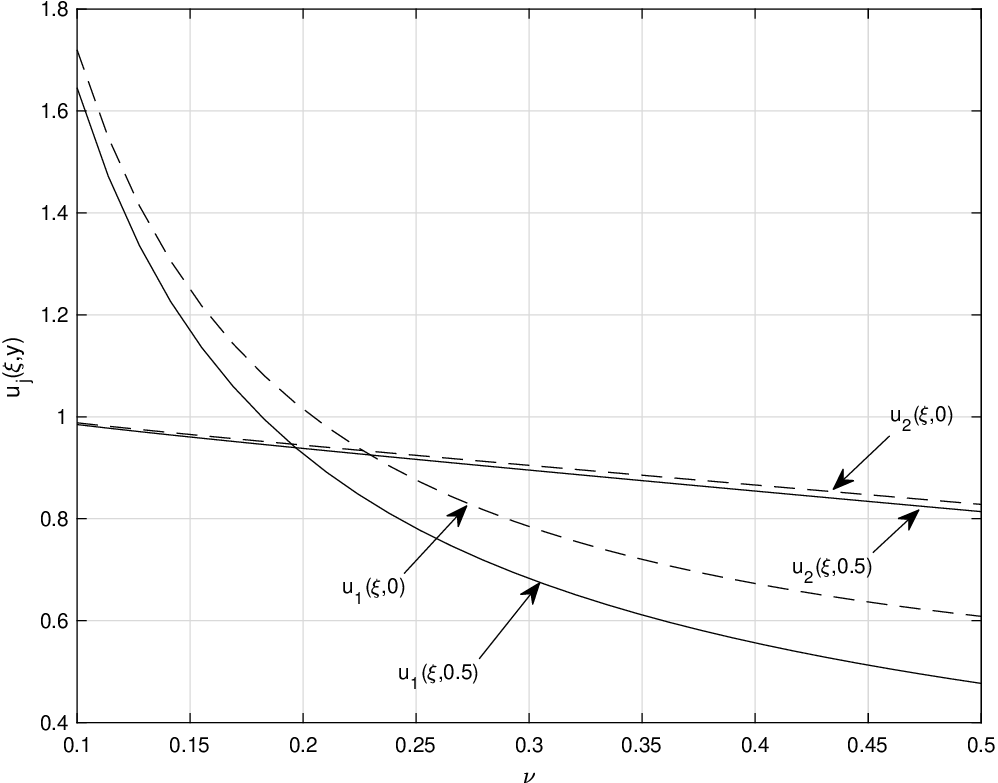}}}
\caption{
The variation of the  tangential and normal displacements  $u_j(\Gx,y)$ for $y=0.5$ (the solid lines) and on the boundary  $y=0$ (dashed lines) when  $0.1\le \nu\le 0.5$,   $\mu_0^{-1}H_1=-1$, $H_2=0$,
$\Gx_0=0$, $\Gx=-1$, $\nu_P=0.3$, $V/c_s=0.2$.}
\label{fig6}
\end{figure} 

\begin{figure}[t]
\centerline{
\scalebox{0.6}{\includegraphics{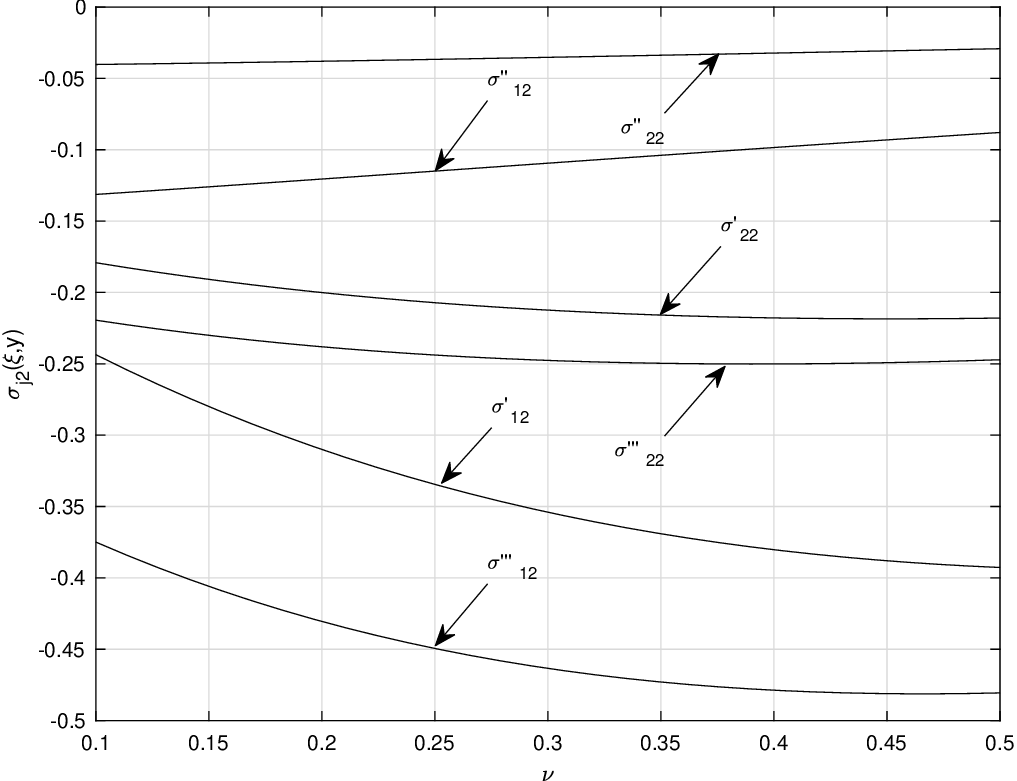}}}
\caption{
The variation of the  tangential and normal stresses $\sigma_{j2}(\Gx,y)$ when $0.1\le \nu\le 0.5$,  $\Gx_0=0$, $\Gx=-1$, $y=0.3$, 
$\nu_P=0.3$, $V/c_s=0.2$ in cases $\mu_0^{-1}H_1=-1$, $H_2=0$ (curves $\Gs'_{12}$ and $\Gs'_{22}$),
$H_1=0$,  $\mu_0^{-1}H_2=-1$  (curves $\Gs''_{12}$ and $\Gs''_{22}$), and $\mu_0^{-1}H_1=-1$, $\mu_0^{-1}H_2=-1$ 
 (curves $\Gs'''_{12}$ and $\Gs'''_{22}$).} 
\label{fig7}
\end{figure} 

\begin{figure}[t]
\centerline{
\scalebox{0.6}{\includegraphics{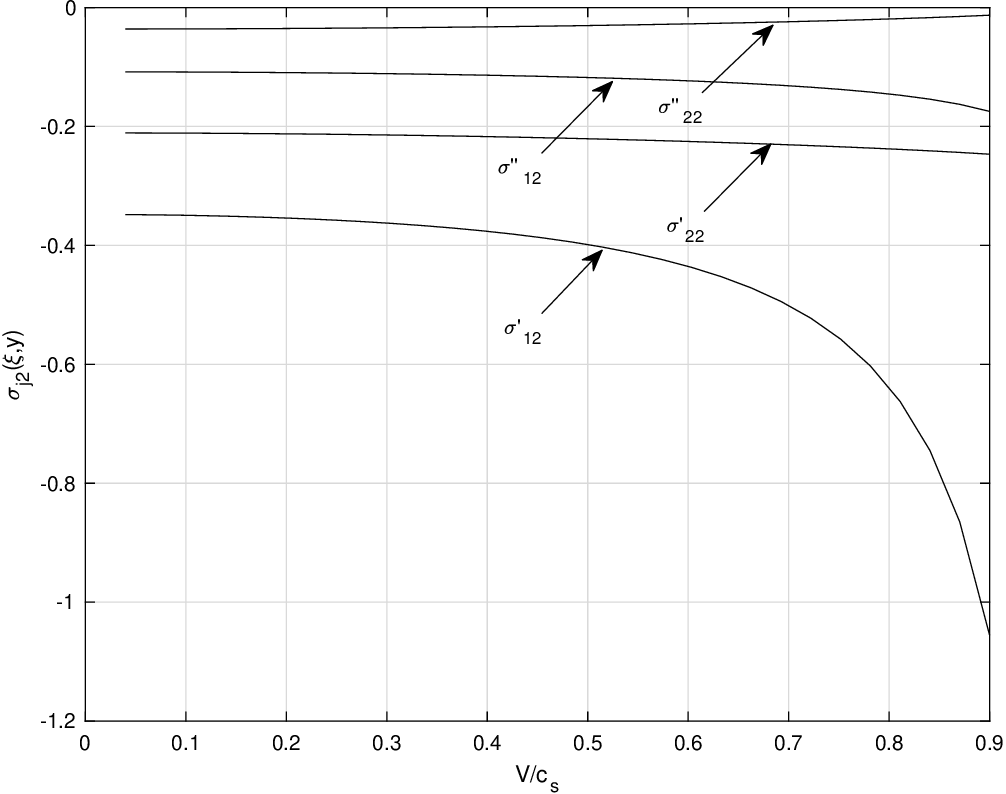}}}
\caption{
The tangential and normal stresses $\sigma_{j2}(\Gx,y)$ versus  the dimensionless speed $V/c_s$ 
in cases $\mu_0^{-1}H_1=-1$, $H_2=0$ (curves $\Gs'_{12}$ and $\Gs'_{22}$) and
$H_1=0$,  $\mu_0^{-1}H_2=-1$  (curves $\Gs''_{12}$ and $\Gs''_{22}$) when $\Gx_0=0$, $\Gx=-1$, $y=0.3$, $\nu=0.3$,
$\nu_P=0.3$.} 
\label{fig8}
\end{figure}

Numerical tests implemented reveal some remarkable properties of the coefficients and the solution of the system of algebraic 
equations  (\ref{6.33}). It turns out that the functions $G_1(s)$ and $G_2(s)$ possess the symmetry property 
\beq
G_j(s)=\ov{G_j(\bar s)}, \quad s=\Gs+it\in\GO.
\label{7.1}
\eeq
Due to the structure of the system of integral equations (\ref{6.4}) in the case $m=1$ we have 
$$
F_{1k+}^{(1\pm)}=F_{1k-}^{(1\pm)}, \quad F_{2k+}^{(1\pm)}=-F_{2k-}^{(1\pm)},
$$
\beq
F_{1k\pm}^{(1+)}=\ov{F_{1k\pm}^{(1-)}}, \quad F_{2k\pm}^{(1+)}=-\ov{F_{2k\pm}^{(1-)}}.
\label{7.2}
\eeq
In the case $m=2$ the analogs of these properties have the form
$$
F_{1k+}^{(2\pm)}=-F_{1k-}^{(2\pm)}, \quad F_{2k+}^{(2\pm)}=F_{2k-}^{(2\pm)},
$$
\beq
F_{1k\pm}^{(2+)}=-\ov{F_{1k\pm}^{(2-)}}, \quad F_{2k\pm}^{(2+)}=\ov{F_{2k\pm}^{(2-)}}.
\label{7.3}
\eeq
It is found that the determinants $\GD_+$ and $\GD_-$ are the same and approximately real. Their numerical values are stable as the number
of equations $4N$ in the system (\ref{6.33}) grows. For $\nu=0.1$, $V/c_s=0.2$, $\nu_P=0.3$, and $\mu_0^{-1}H_1=\mu_0^{-1}H_2=-1$ we have
$$
\GD_\pm=0.9821-i2.013\times 10^{-4} \;\; {\rm for}\;  N=50,
$$
$$
\GD_\pm=0.9944-i1.6671 \times 10^{-4} \;\; {\rm for}\;  N=75,
$$
\beq
\GD_\pm=1.0007-i1.4402  \times 10^{-4} \;\; {\rm for}\;  N=100.
\label{7.4}
\eeq
With the same level of accuracy the coefficients $C_{1\pm}$ and $C_{2\pm}$ have the properties
\beq
C_{1+}=\ov{C_{1-}}, \quad C_{2+}=\ov{C_{2-}},
\label{7.5}
\eeq
and the values of the displacements and the stresses are real. The results have to be invariants of the parameter $\Gs$ provided it
falls in the interval $(0,\nu)$. This is confirmed by numerical tests: for $N=100$ and the data used in recovering the determinants 
$\GD_\pm$ we
have $\GD_\pm=1.0007-i1.4402  \times 10^{-4}$ if $\Gs=\nu/4$ and $0.9953-i1.5088 \times 10^{-4}$ if $\Gs=\nu/2$.

To reconstruct the displacements, their derivatives, and the stresses, we use the asymptotic expansions (\ref{5.17'}) and (\ref{5.25}).
The coefficients $d_{j0}$ and $d_{j2}$ in these expansions need the values of the constants $C_{j\pm}$ and therefore the values of
$\GF_{j\pm}^{(m)}(\nu-1)$, $j,m=1,2$, given by
\beq
\GF_{j\pm}^{(m)}(\nu-1)=\pm\fr{1}{4}\int_\GO\fr{G_j(p)\GF^{(m)}_{3-j\pm}(p-1)dp}{\cos\fr{\pi}{2}(p-\nu)}
-\fr{\Gd_{jm}}{\cos\fr{\pi\nu}{2}}, \quad j,m=1,2.
\label{7.6}
\eeq
By the method applied in Section 6 we evaluate the integral approximately and express it through the solution of the algebraic system
 (\ref{6.33}). We have 
$$
\GF_{j\pm}^{(m)}(\nu-1)=\pm\fr{i}{2\pi}\sum_{n=1}^N
\left\{\fr{F_{3-j\, n\pm}^{(m+)}\left(x_n^{i\Gd^+_{3-j}+1}-x_{n-1}^{i\Gd^+_{3-j}+1}\right)}
{[x_n e^{-i\pi(\Gs-\nu)/2}+e^{i\pi(\Gs-\nu)/2}](i\Gd_{3-j}^+ +1)}
\right.
$$
\beq
\left.
+\fr{F_{3-j\, n\pm}^{(m-)}\left(x_n^{i\Gd^-_{3-j}+1}-x_{n-1}^{i\Gd^-_{3-j}+1}\right)}
{[x_n e^{i\pi(\Gs-\nu)/2}+e^{-i\pi(\Gs-\nu)/2}](i\Gd_{3-j}^- +1)}\right\}-\fr{\Gd_{jm}}{\cos\fr{\pi\nu}{2}}.
\label{7.7}
\eeq

For all numerical tests we take the Poisson ratio $\nu_P=0.3$, $\Gx=-1$, and $\Gx_0=0$.
It appears that the displacements on the boundary of the half-plane are decreasing as the parameter $\nu$ 
is growing while the concentrated force runs at constant speed. The two displacements are plotted in Fig.1 
for $V=0.2c_s$, $H_1/\mu_0=H_2/\mu_0=-1$. If we replot these curves for the tangential derivatives of the 
displacements (Fig. 2), we find out that the tendency is reverse: the derivatives grow when  the parameter $\nu$ is growing.
In Figs. 3 and 4 for the values $0.2$, $0.3$, and $0.4$ of the parameter $\nu$ and the same concentrated load, $h_j(\Gx)=H_j\Gd(\Gx-\Gx_0)$ with 
$H_1/\mu_0=H_2/\mu_0=-1$, we plot the displacements $u_1(\Gx,0)$ and $u_2(\xi,0)$ as functions of  speed varying  in the subsonic range. It is seen that both displacements grow as the speed grows. In Fig. 5 we compare the displacements as functions of the parameter $\nu$ on the boundary with
their values at the point  $\xi=-1$, $y=0.5$. These values are recovered  by using the asymptotic expansion (\ref{5.17'}). The speed $V=0.2 c_s$
and the concentrated load applied at $\Gx_0=0$ is characterized by  $H_1/\mu_0=-1$ and $H_2=0$. 

In Fig. 6 we show the results of computations of the stresses $\Gs_{12}$ and $\Gs_{22}$ as the point $\xi=-1$, $y=0.3$.
As before, we take  $V=0.2 c_s$, $\xi_0=0$, $\nu\in[0.1,0.5]$. The curves are presented for three cases, $H_1/\mu_0=-1$, $H_2=0$,
$H_1=0$, $H_2/\mu_0=-1$, and $H_1/\mu_0=H_2/\mu_0=-1$.
In Fig. 7 in two cases, $H_1/\mu_0=-1$, $H_2=0$ and
$H_1=0$  $H_2/\mu_0=-1$, and when $\nu=0.3$  we plot the stresses at the same point as functions of the dimensionless speed $V/c_s$.

  \section{Conclusions} 

We have solved a steady state two-dimensional model problem of an inhomogeneous plane subjected to a load running
along the boundary at subsonic speed when the Lame coefficients and the density are power functions of depth. 
The methodology applied is based on the Fourier and Mellin transform and analysis of the resulting Carleman boundary
value problem for two meromorphic functions in a strip with two shifts. We have managed to express the unknown functions
through the solution of a system of four singular integral equations on the interval $(0,1)$ with a fixed singularity and oscillating coefficients.
We have developed a numerical method for its solution. Numerical tests have demonstrated its efficiency and a good accuracy.

There are several differences between the model proposed in this paper and the one traditionally used in contact mechanics
to describe the contact interaction of a stamp and a power-law graded foundation.  
Firstly, the classical model recovers the pressure distribution
only in the contact zone and is not capable to predict the stress and displacement fields in the interior of the 
power-law graded half-plane.   In the framework of the new model we have determined integral representations
for the mechanical fields everywhere in the half-plane and derived asymptotic expansions  convenient near the boundary and far away from it.
Secondly, the model proposed is steady state and, by passing to the limit $V\to 0$, gives the solution to the static problem.
It also admits a generalization to the transient case, while the traditional model is static and is not applicable for the dynamic case.  Lastly, the standard model employs the Flamant model solution 
$\Gs_r=C r^{-1}\GF$, where $C$ is a constant, $\GF$ is a function of $\Gvf$, and $(r,\Gvf)$ are polar coordinates.
The {\it ab initio} model proposed is based on the
momentum balance equations of dynamic elasticity,  the stress-strain relations, and  the traction boundary conditions
$\Gs_{j2}(\Gx)=h_j(\Gx)$, $\Gx=x-Vt$. Due to the Lame coefficients representations $\Gl(y)=\Gl_0y^\nu$,
$\Gm(y)=\Gm_0y^\nu$, $0<\nu<1$,  the relations (\ref{2.8}) between  the stresses and displacements, and also because 
the displacements are bounded on the surface $y=0$, the  boundary conditions (\ref{2.9}) are not affected by the tangential derivatives 
$\Md u_j/\Md\Gx$,
while the strains $\Gve_{ij}$ have a power singularity of order $\nu$ as $y\to 0^+$. This implies that, when 
$\nu\to 0^+$,  the boundary conditions (\ref{2.1}) do not tend to the boundary conditions
$$
\mu_0\left(\fr{\Md u_1}{\Md y}+\fr{\Md u_2}{\Md\Gx}\right)=h_1(\Gx),
\quad 
\Gl_0\fr{\Md u_1}{\Md \Gx}+(\Gl_0+2\mu_0)\fr{\Md u_2}{\Md y}=h_2(\Gx),\quad |\Gx|<\infty,
$$
of the homogeneous case. That is why the solution in the case $\nu>0$ does not tend to the homogeneous solution  
as $\nu\to 0^+$. 

The model problem of a load running along the boundary of a power-law graded half-plane solved in this work
 has a potential to be used as a Green function in modeling of static and dynamic contact 
interaction of a stamp with a power-law graded foundation and also in studying crack propagation along the interface between two 
power-law graded materials. It is of interest to continue research in this direction.

 \vspace{.2in}

{\centerline{\Large\bf  References}}

\vspace{.1in}

\begin{enumerate}

 \item\label{lek} S. G. Lekhnitskii,  Radial distribution of stresses in a wedge and in a half-plane with variable modulus of elasticity.
\textit{J. Appl. Math. Mech. (PMM)}   \textbf{26} 1962,  199-206.

 \item\label{ros} N. A. Rostovtsev,  On the theory of elasticity of a nonhomogeneous medium. \textit{J. Appl. Math. Mech. (PMM)}  \textbf{28} 1964, 745-757.

 \item\label{kor}  B. G. Korenev,  A die resting on an elastic half-space, the modulus of elasticity of which is an
exponential function of depth.  \textit{Dokl. Akad. Nauk SSSR}  \textbf{112} 1957, 823-826.

 \item\label{mos}  V. I. Mossakovskii, Pressure of a circular die  on an elastic half-space, whose modulus of elasticity is an exponential  function of depth.   \textit{J. Appl. Math. Mech. 
(PMM)}   \textbf{22} 1958,  168-171. 

 \item\label{pop1} G. Ia. Popov, On a method of solution of the axisymmetric contact problem of the theory of elasticity.
 \textit{J. Appl. Math. Mech. (PMM)}   \textbf{25}  1961, 105-118.

 \item\label{pop2} G. Ya. Popov,  \textit{Contact problems for a linearly deformed base}, Vishcha Shkola, Kiev,1982.

\item\label{joh}  K. L. Johnson,  K. Kendall, A. D. Roberts,   Surface energy and the contact of elastic solids. 
 \textit{Proc Roy. Soc. A } \textbf{324} 1971, 301-313.

 \item\label{gia} A. E. Giannakopoulos, P. Pallot,   Two-dimensional contact analysis of elastic
graded materials. \textit{J. Mech. Phys. Solids} \textbf{48} 2000, 1597-1631.

 \item\label{che2} S. Chen, C. Yan, P.  Zhang, H. Gao,  
 Mechanics of adhesive contact on a power-law graded elastic half-space. \textit{J. Mech. Phys. Solids} \textbf{57} 2009, 1437-1448.

 \item\label{guo} X. Guo, F. Jin, H.  Gao,   Mechanics of non-slipping adhesive contact on a power-law graded elastic
half-space.  \textit{Int.  J. Solids Struct.} \textbf{48} 2011, 2565-2575.

 \item\label{hes} M. Hess, A simple method for solving adhesive and non-adhesive axisymmetric contact problems of elastically graded materials. \textit{Int.  J.  Eng. Sci.}  \textbf{104} 2016, 20-33.

 \item\label{ant1} Y. A. Antipov,  S. M. Mkhitaryan,  Axisymmetric contact of two different power-loaw gradede elastic bodies and an integral equation with two Weber-Schafheitlin kernels.  \textit{Quart. J. Mech. Appl. Math.}  \textbf{75} 2022, 393-420. 

 \item\label{ant2} Y. A. Antipov,  S. M. Mkhitaryan,  Hertzian and adhesive plane models of contact of two inhomogeneous elastic bodies.
 \textit{Eur. J. Appl. Math.}  \textbf{34} 2023, 667-700.

 \item\label{ant3}  Y. A. Antipov, Plane and axisymmetric problems of an interfacial crack in a power-law graded material.
 \textit{Proc. Roy. Soc.}  \textbf{479} (2279) 2023, Paper No. 20230630. 
 
 \item\label{car}  T. Carleman  Sur la th\'eorie des \'equations int\'egrales et ses applications. Verhandl. des
Internat. Mathem. Kongr., 1932, 138–151.

 \item\label{zve} E. I. Zverovich, G. S.  Litvinchuk,
Boundary value problems with shift for analytic functions, and singular functional equations.
 \textit{Uspehi Mat. Nauk}  \textbf{23} 1968, no.3, 67-121.
 
 \item\label{che} Y. I. \v Cerski\u i,
A normally solvable equation of smooth transition.
 \textit{Dokl. Akad. Nauk SSSR} \textbf{190} 1970, 57-60.
 
 \item\label{ban}  R. D. Bancuri, 
A contact problem for a wedge with elastic bracing.
 \textit{Dokl. Akad. Nauk SSSR} \textbf{211} 1973, 797-800.
 
 \item\label{ant4}  Y. A. Antipov and V. V. Silvestrov, Second-order functional-difference equations. I.: Method of the Riemann-Hilbert problem on Riemann surfaces. \textit{Quart. J. Mech. Appl. Math.}  \textbf{57} 2004, 245-265. 
 
 \item\label{ant5}  Y. A. Antipov and V. V. Silvestrov, Vector functional-difference equation in electromagnetic scattering.  \textit{IMA J. Appl. Math.}, \textbf{69} 2004, 27-69. 

 \item\label{ant6} Y. A. Antipov, Diffraction of a plane wave by a circular cone with an impedance boundary condition.
  \textit{SIAM J. Appl. Math.}  \textbf{62} 2002,  1122-1152. 

 \item\label{ant7} Y. A. Antipov,  V. V. Silvestrov, Method of integral equations for systems of difference equation.   \textit{IMA 
 J. Appl. Math.}   \textbf{72} 2007, 681-705. 

\item\label{gra} I. S.  Gradshteyn, I.M.  Ryzhik,  \textit{Table of Integrals, Series, and Products}, Academic Press, New York, 1994.

 \end{enumerate}

 \end{document}